\documentclass[11pt, letterpaper, oneside]{article}

\usepackage{latexsym, amsmath, amssymb, amsfonts, amscd}
\usepackage{amsthm}
\usepackage{t1enc}
\usepackage[mathscr]{eucal}
\usepackage{indentfirst}
\usepackage{graphicx}
\usepackage{fancyhdr}
\usepackage{fancybox}
\usepackage{enumerate,psfrag}
\usepackage[all,poly,web,knot]{xy}
\usepackage{epsfig,diagrams}

\theoremstyle{plain}
\newtheorem{thm}{Theorem}[section]

\newtheorem{prop}[thm]{Proposition}
\newtheorem{lemma}[thm]{Lemma}
\newtheorem{cor}[thm]{Corollary}

\newtheorem{claim}[thm]{Claim}

\renewcommand{\latticebody}{\drop@{ }}

\theoremstyle{definition}
\newtheorem{defi}[thm]{Definition}
\newtheorem{pdef}[thm]{Proposition-Definition}

\theoremstyle{remark}
\newtheorem{remark}[thm]{Remark}

\newcommand{\la}{\leftarrow}
\newcommand{\ra}{\rightarrow}
\newcommand{\lra}{\longrightarrow}
\newcommand{\lla}{\longleftarrow}
\newcommand{\rra}{\Rightarrow}

\newcommand{\thla}{\twoheadleftarrow}
\newcommand{\thra}{\twoheadrightarrow}

\newcommand{\N}{\ensuremath{\mathbb N}}

\newcommand{\R}{\ensuremath{\mathbb R}}

\renewcommand{\d}{d}

\newcommand{\D}{\mathscr D}

\newcommand{\cA}{\mathcal{A}}

\newcommand{\cX}{\mathcal{X}}
\newcommand{\cY}{\mathcal{Y}}

\newcommand{\cG}{\mathcal{G}}
\newcommand{\cH}{\mathcal{H}}

\newcommand{\bepsilon}{\mbox{\boldmath $\epsilon$}}
\newcommand{\bDelta}{\mbox {\boldmath $\Delta$}}

\newcommand{\tG}{\tilde{G}}

\newcommand{\tU}{\tilde{U}}

\newcommand{\teta}{\tilde{\eta}}

\newcommand{\be}{\bar{e}}
\newcommand{\bg}{\bar{g}}
\newcommand{\bareta}{\bar{\eta}}
\newcommand{\bgamma}{\bar{\gamma}}

\newcommand{\bt}{\mathbf{t}}                  
\newcommand{\bs}{\mathbf{s}}                  
\newcommand{\bbt}{\bar{\mathbf{t}}}           
\newcommand{\bbs}{\bar{\mathbf{s}}}           

\def\R{{\mathbb R}}

\def\N{{\mathbb N}}
\def\U{{\mathcal U}}
\def\L{\Lambda}
\def\D{\Delta}
\newarrow{into} C--->
\newarrow{Eq} =====
\newarrow{dashto} ....>

\def\pD{\partial\D}
\def\lht{Lie homotopy type}

\def\PB(#1,#2,#3,#4){
\left\{\begin{matrix}#1&\!\!\!\stackrel{?}{\longrightarrow}&\!\!\!#2\\
\downarrow&&\!\!\!\downarrow\\
#3&\!\!\!\stackrel{?}{\longrightarrow}&\!\!\!#4\end{matrix}\right\}}

\newtheorem{df}[thm]{Definition}

\begin{document}

\title{Lie $n$-groupoids and stacky Lie groupoids}
\author{Chenchang Zhu \thanks{Research partially supported by the Liftoff fellowship 2004 of the Clay
Institute}\\
 Departement Mathematik\\
    ETH Zentrum, R\"amistrasse 101\\
    8092 Z\"urich, Switzerland \\
  \small{(zhu@math.ethz.ch) }}
\date{\today}

\maketitle

\begin{abstract}
We discuss two generalizations of Lie groupoids. One consists of
Lie $n$-groupoids defined as simplicial manifolds with trivial
$\pi_{k\geq n+1}$. The other consists of stacky Lie groupoids
$\cG\rra M$ with $\cG$ a differentiable stack. We build a 1-1
correspondence between Lie 2-groupoids and stacky Lie groupoids up
to a certain Morita equivalence. Equivalences of these higher
groupoids are also described.
\end{abstract}

\tableofcontents

\section{Introduction}

Recently there has been much interest in higher group(oid)s, which
generalize the notion of group(oid)s in various ways. Some of them
turn out to be unavoidable to study problems in differential
geometry. An example comes from the string group, which is a
3-connected cover of $Spin(n)$. More generally, to any compact
simply connected group $G$ one can associate its string group
$String_G$. It has various models, given by Solz and Teichner
\cite{stolz} \cite{st} using an infinite dimensional extension of
$G$, by Brylinski \cite{bry-mc1} using a $U(1)$-gerbe with the
connection over $G$, and recently by Baez et al.
\cite{baez:str-gp} using Lie 2-groups and Lie 2-algebras.
Henriques \cite{henriques} constructs the string group as a higher
group in the setting of this paper and as an integration object of
a certain Lie 2-algebra using an integration procedure similar to
that of \cite{s-funny, z2}.

Other examples come from  a kind of \'etale stacky groupoids
(Weinstein groupoids) \cite{tz}. They are the 1-1 global objects
integrating Lie algebroids, where a Lie algebroid could be roughly
understood as a bundle of Lie algebras. Notice that unlike (finite
dimensional) Lie algebras which always have associated Lie groups,
Lie algebroids do not always have  associated Lie groupoids
\cite{am1, am2}. One needs to enter the world of stacky groupoids
to obtain the desired 1-1 correspondence. Since Lie algebroids are
closely related to Poisson geometry, this result applies to
complete the first step of Weinstein's program of quantization of
Poisson manifolds: to associate to Poisson manifolds their
symplectic groupoids \cite{w-poisson, wx}. It turns out that some
``non-integrable'' Poisson manifolds cannot have symplectic (Lie)
groupoids. This problem is solved in \cite{tz2} with the above
result so that every Poisson manifold has a corresponding \'etale
stacky symplectic groupoid.

Higher group(oid)s were already studied in the early twentieth
century by Whitehead and his followers under various terms, such
as cross-modules. However in this paper we use a uniform method to
describe higher (respectively Lie) groupoids  using simplicial
sets \cite{may} (respectively manifolds) since it is believed that
there should be an equivalence between $n$-groupoids and spaces
whose homotopy groups are trivial above $\pi_n$ (also called
$n$-coskeleta or $n$-truncated homotopy types)\footnote{This was suggested
(indirectly) by Jacob Lurie to the author, however it was known much
earlier, for example,  by Duskin
and Glenn \cite{duskin, glenn}.}. 
The 0-simplices of the simplicial set
correspond to the objects, the 1-simplices to the arrows (or
1-morphisms), and the higher dimensional simplices the higher
morphisms. This method  becomes more suitable especially when we
are dealing with the differentiable category.

Recall that a simplicial set (respectively manifold) $X$ is made
up by sets (respectively manifolds) $X_n$ and structure maps
\[ d^n_i: X_n \to X_{n-1} \;\text{(face maps)}\quad s^n_i: X_n \to X_{n+1} \; \text{(degeneracy maps)},\;\; \text{for}\;i\in \{0, 1, 2,\dots, n\} \]
that satisfy the coherence conditions
\begin{equation}\label{eq:face-degen}
\begin{split}
 d^{n-1}_i d^{n}_j = d^{n-1}_{j-1} d^n_i \; \text{if}\; i<j, &\quad s^{n}_i s^{n-1}_j = s^{n}_{j+1} s^{n-1}_i \; \text{if}\; i\leq j,
 \\
 d^n_i s^{n-1}_j = s^{n-2}_{j-1} d^{n-1}_i \; \text{if}\; i<j, &\quad
 d^n_j s^{n-1}_j=\mathrm{id}=d^n_{j+1} s^{n-1}_j,\quad
 d^n_i s^{n-1}_j = s^{n-2}_j d^{n-1}_{i-1} \; \text{if}\; i> j+1.
\end{split}
\end{equation}

The first two examples are the simplicial $m$-simplex $\D[m]$ and
the horn $\L[m,j]$ with
\begin{equation}\label{eq:simplex-horn}
\begin{split}
(\D[m])_n & = \{ f: (0,1,\dots,n) \to (0,1,\dots, m)| f(i)\leq
f(j),
\forall i \leq j\}, \\
(\L[m,j])_n & = \{ f\in (\D[m])_n| \{0,\dots,j-1,j+1,\dots,m\}
\nsubseteq \{ f(0),\dots, f(n)\} \}.
\end{split}
\end{equation}
In fact the horn $\L[m,j]$ is a simplicial set obtained from the
simplicial $m$-simplex $\D[m]$ by taking away its unique
non-degenerate $m$-simplex as well as the $j$-th of its $m+1$
non-degenerate $(m-1)$-simplices, as in the following picture (in
this paper all the arrows are oriented from bigger numbers to
smaller numbers): \vspace{.6cm}

\centerline{\epsfig{file=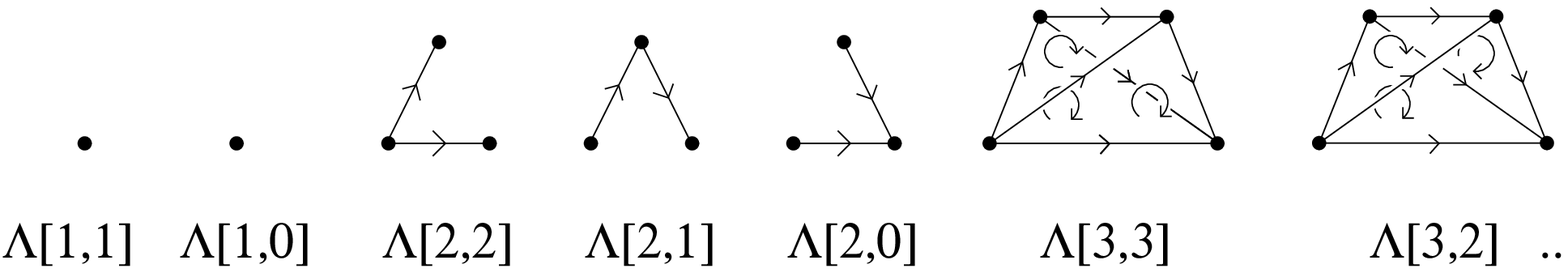,height=1.9cm}}
\vspace{.6cm}

A simplicial set $X$ is {\em Kan} if any map from the horn
$\L[m,j]$ to $X$ ($m\ge 1$, $j=1,\dots,m$), extends to a map from
$\D[m]$. Let us call $Kan(m,j)$ the Kan condition for the horn
$\L[m,j]$. A Kan simplicial set is therefore a simplicial set
satisfying $Kan(m,j)$ for all $m\ge 1$ and $0\leq j\leq m$. In the
language of groupoids, the Kan condition corresponds to the
possibility of composing various morphisms. For example, the
existence of a composition for arrows is given by the condition
$Kan(2,1)$, whereas the composition of an arrow with the inverse
of another is given by $Kan(2,0)$ and $Kan(2,2)$. \vspace{.6cm}

\begin{equation}\label{compo}
\epsfig{file=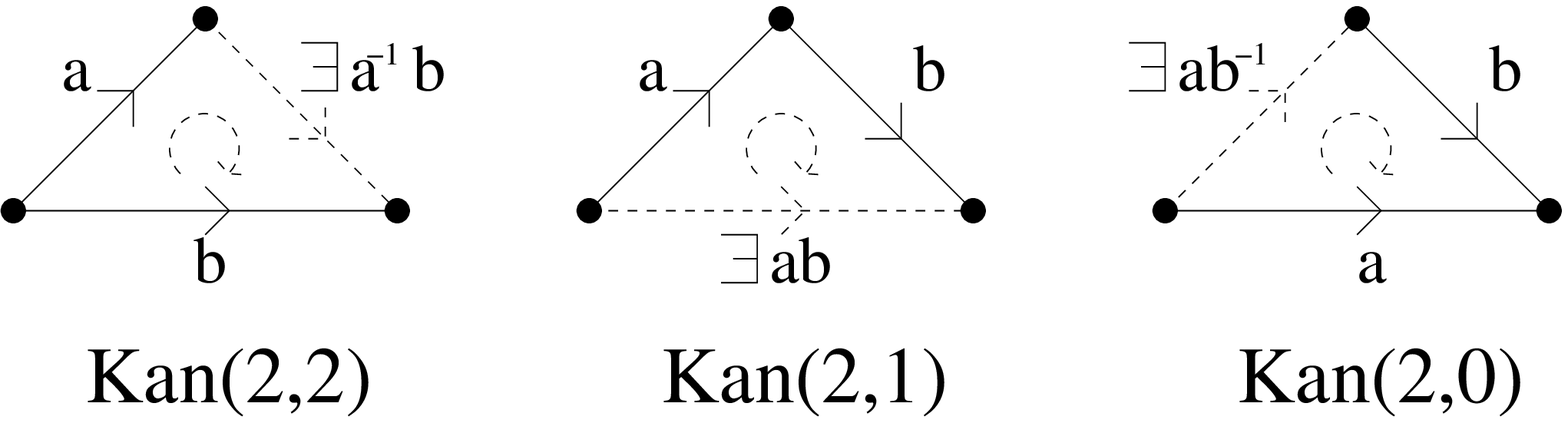,height=2cm}
\end{equation}

Note that the composition of two arrows is in general not unique,
but any two of them can be joined by a 2-morphism $h$ given by
$Kan(3,1)$. \vspace{.6cm}

\begin{equation}\label{composition}
\epsfig{file=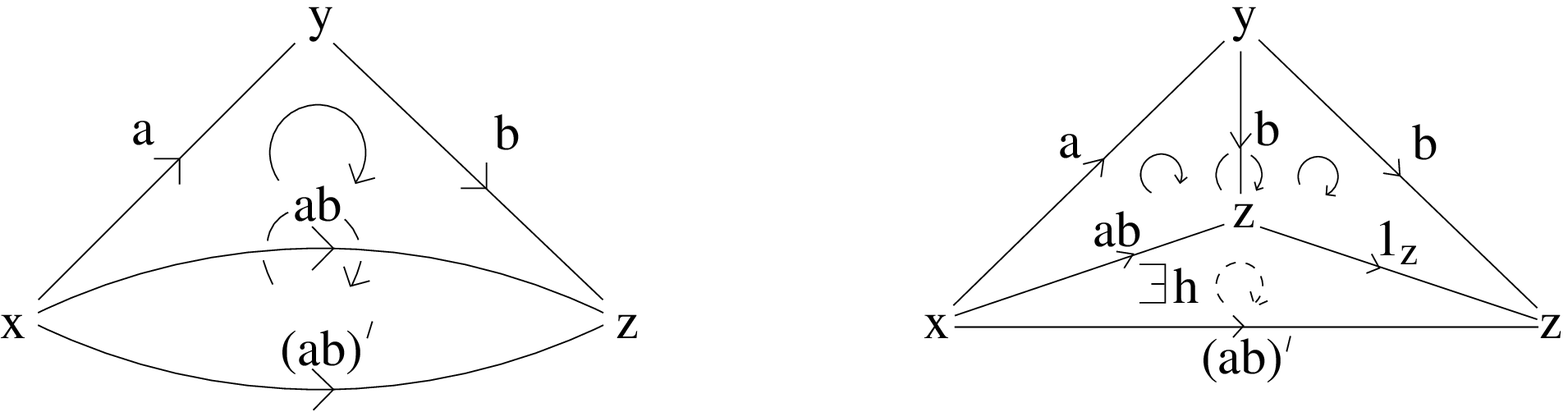,height=2.5cm}
\end{equation}

Here, $h$ ought to be a bigon, but since we don't have any bigons
in a simplicial set, we view it as a triangle with one of its
edges degenerate. The degenerate 1-simplex above $z$ is denoted
$1_z$.

In an $n$-groupoid, the only well defined composition law is the
one for $n$-morphisms. This motivates the following definition.

\begin{df}\label{groupoid}
An $n$-groupoid ($n\in\N \cup \infty$) $X$ is a simplicial set
that satisfies $Kan(m,j)$ for all $0\le j\le m\ge 1$ and
$Kan!(m,j)$ for all $0\le j\le m<n$, where\\
\noindent
\begin{tabular}{p{1.6cm}p{10cm}}
$Kan(m,j)$:&  Any map $\L[m,j]\to X$ extends to a map $\D[m]\to X$.\\
$Kan!(m,j)$: &
Any map $\L[m,j]\to X$ extends to a unique map $\D[m]\to X$.\\
\end{tabular}\\
An $\infty$-groupoid will be called a homotopy type.  
\end{df}

An {\em $n$-group} is an $n$-groupoid for which  $X_0$ is a point.
When $n=2$, they are the same as the weak $2$-group(oid)s in
\cite{noohi} but different from the various kinds of
$2$-group(oid)s or double groupoids in \cite{bala:2gp,
br-sp}\footnote{ See \cite{henriques} for an explanation of the
relation with \cite{bala:2gp}.}. A usual groupoid (category with
only isomorphisms) is equivalent to a 1-groupoid in the sense of
Definition \ref{groupoid}. Indeed, from a usual groupoid, one can
form a simplicial set whose $n$-simplices are given by sequences
of $n$ composable arrows. This is a standard construction called
the \em nerve \rm of a groupoid and one can check that it
satisfies the required Kan conditions.

On the other hand, a 1-groupoid $X$ in the sense of Definition
\ref{groupoid} gives us a usual groupoid with objects and arrows
given respectively by the 0-simplices and 1-simplices of $X$. The
unit is provided by the degeneracy $X_0\to X_1$, the inverse and
composition are given by the Kan conditions $Kan(2,0)$, $Kan(2,1)$
and $Kan(2,2)$ as in (\ref{compo}), and the associativity is given
by $Kan(3,1)$ and $Kan!(2,1)$. \vspace{.6cm}

\centerline{\epsfig{file=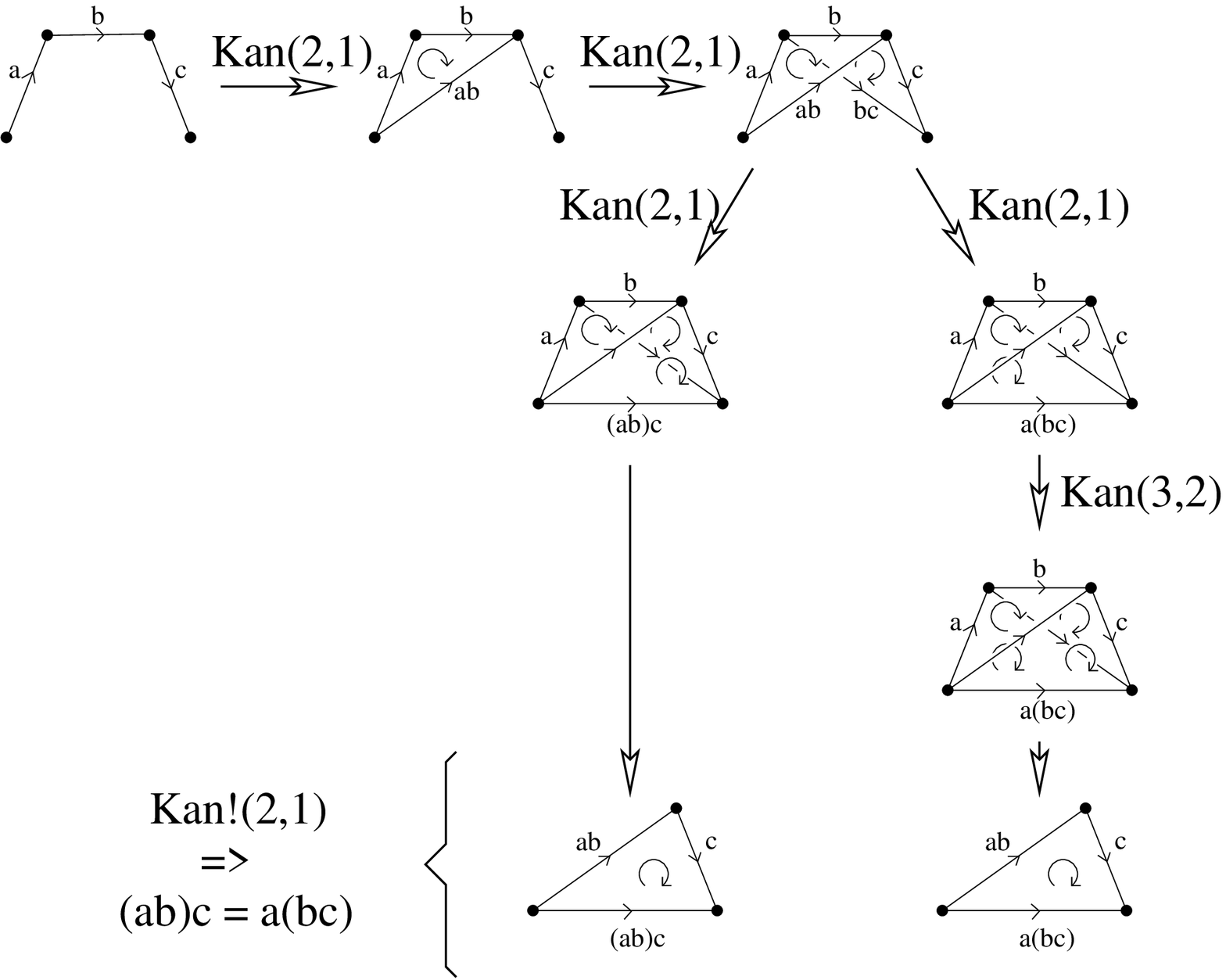,height=9cm}}
\centerline{Proof of associativity.} \vspace{.6cm}

This motivates the corresponding definition in the differentiable
category.

\begin{df}\label{defngroupoids}
A Lie $n$-groupoid $X$ ($n\in\N \cup \infty$) is a simplicial
manifold that satisfies $Kan(m,j)$ for all $0\le j\le m\ge 1$ and
$Kan!(m,j)$ $0\le j\le m<n$, where\\ \noindent
\begin{tabular}{p{1.6cm}p{9.6cm}}
$Kan(m,j)$:& The restriction map $\hom(\D[m],X)\to\hom(\L[m,j],X)$
is a
surjective submersion.\\
$Kan!(m,j)$:& The restriction map
$\hom(\D[m],X)\to\hom(\L[m,j],X)$ is a
diffeomorphism.\\
\end{tabular} \\A Lie $\infty$-groupoid will be
called a Lie homotopy type. \vspace{.3cm}

\end{df}

We view simplicial sets $\D[m]$ and $\L[m,j]$ as simplicial
manifolds with their discrete topology. Then $\hom(S,X)$ denotes
the set of homomorphisms of simplicial manifolds, with its natural
topology, so $\hom(\D[m],X)$ is just another name for $X_m$.
However the fact that $\hom(\L[m,j],X)$ is a manifold is not
obvious (see Section \ref{sect:lie-n-gpd}).

On the other hand, a {\em stacky Lie (SLie) groupoid} $\cG \rra
M$, following the concept of Weinstein (W-) groupoid in \cite{tz},
is a groupoid object in the world of differentiable stacks with
its base $M$ an honest manifold (see Definition \ref{def:sliegpd}
for detailed axioms). When $\cG$ is also a manifold, $\cG \rra M$
is obviously a Lie groupoid.  {\em W-groupoids}, which are {\em
\'etale SLie} groupoids, provide a way to build the 1-1
correspondence with Lie algebroids.

Given these two  higher generalizations of Lie groupoids, Lie
$n$-groupoids and SLie groupoids, arising from different
motivations and constructions, we ask the following questions:
\begin{itemize}
\item Are SLie groupoids the same as Lie $n$-groupoids for some $n$?
\item If not exactly, to which extent they are the same?
\item Is there a way to also realize Lie $n$-groupoids as
integration objects of Lie algebroids?
\end{itemize}

In this paper, we answer the two first questions by

\begin{thm} \label{2-w}
There is a one-to-one correspondence between SLie (respectively
W-) groupoids and Lie 2-groupoids (respectively Lie 2-groupoids
whose $X_2$ is \'etale over $\hom(\L[2,j],X)$) modulo 1-Morita
equivalences\footnote{Morita equivalences preserving $X_0$} of Lie
2-groupoids.
\end{thm}

The last question will be answered positively in a future work \cite{z2}:

\begin{thm}\label{thm:2-a}
Let $A$ be a Lie algebroid and let $Lmor(-,-)$ be the space of Lie
algebroid homomorphisms satisfying suitable boundary conditions.
Then
\[ Lmor(T\Delta^2, A)/ Lmor(T\D^3 , A) \Rrightarrow Lmor(T\Delta^1, A) \rra
Lmor(T\Delta^0, A), \] is a Lie 2-groupoid corresponding to the
W-groupoid $\cG(A)$ constructed in \cite{tz} under the
correspondence in the above theorem.
\end{thm}

For the first theorem, we (are forced to) develop Morita
equivalence of Lie 2-groupoids, which  is expected to be useful in
the theory of 2-stacks and 2-gerbes and should correspond to
Morita equivalence of SLie groupoids \cite{bz}.

We also find some technical improvements of the concept of SLie
groupoids: it turns out that an SLie groupoid $\cG\rra M$ always
has a ``good groupoid presentation'' $G$ of $\cG$, which possesses
a strict groupoid map $M\to G$. Moreover the condition on the
inverse map could be simplified.

\noindent {\bf Acknowledgements:} Here I would like to thank
Henrique Bursztyn, Ezra Getzler, Alan Weinstein,
for their hosts and very helpful discussions. I also thank Laurent
Bartholdi and Marco Zambon for many editting suggestions.  Finally I thank especially Andr\'e Henriques who pointed out to me the
potential correspondence of stacky groupoids and Lie 2-groupoids
during the conference of ``Groupoids and Stacks in Physics and
Geometry'' in CIRM-Luminy 2004. I owe a lot to discussions with
Andr\'e. He contributed the exact definitions of  Lie n-groupoids
and their equivalences,  and also nice pictures!

\section{Lie $n$-groupoids} \label{sect:lie-n-gpd}

In differential geometry, Lie groupoids have been studied a lot
(See \cite{cw} for details). They are used to study foliations,
and more recently orbifolds and differentiable stacks
\cite{m-orbi} \cite{bx1}. Here we will try to convince the reader
that it is fruitful to consider them within the context of Lie
$n$-groupoids (Def. \ref{defngroupoids}), especially if one wants
to define and use sheaf cohomology. We shall also embed the
category of manifolds in the category of Lie homotopy types by
sending a manifold $M$ to the constant simplicial manifold
$\underline{M}$ with $\underline{M}_n=M$. Note that the constant
simplicial sets (the ones whose only non-degenerate simplices are
0-dimensional) are exactly the 0-groupoids and they form a
category which is equivalent to the category of sets. Similarly,
the constant simplicial manifolds are exactly the Lie 0-groupoids.
In the future, we shall abuse language and say ``$X$ is a
manifold'' to mean ``X is a constant simplicial manifold''.

As promised after Def. \ref{defngroupoids}, we first prove the
important fact that $\hom(\L[j,m],X)$ is a manifold to make the
definition complete. Since Lie $n$-groupoids are special Lie
homotopy types (or you could call them Kan simplicial manifolds),
we only have to show this fact for Lie homotopy types.

\begin{lemma}\label{collapsable1}
Let $S$ be a finite $k$-dimensional collapsible simplicial set,
namely one that can be obtained from the point by successively
filling horns. Let $X$ be a simplicial manifold that satisfies
$Kan(m,j)$, $\forall 1\le m \le k$, $j=0,\dots,m$. Then the space
$\hom(S,X)$ is naturally a manifold.
\end{lemma}

Here we assume that the spaces $\hom(\L[m,j],X)$ are known to be
manifolds when $m\le k$, for the condition $Kan(m,j)$ to make
sense.

\begin{proof} Let $S'\subset S$ be a collapsible sub-simplicial set such that
$S$ is obtained from $S'$ by filling one horn. In other words,
there is a push-out diagram
\begin{diagram}
S'      &    \rTo    & \SWpbk S\\
 \uTo       &   & \uTo \\
\L[m,j]&    \rinto    &\D[m]&,
\end{diagram}
where the lower arrow is the inclusion. When applying
$\hom(\;\;,X)$ to the above push-out, we get a pull-back diagram

\begin{diagram}
\hom(S',X)&     \lTo    & \hom\SWpbk(S,X)\\
 \dTo       &    & \dTo \\
\hom(\L[m,j],X)& \lTo& \hom(\D[m],X)\;&=X_m
\end{diagram}

By induction on the number of simplices of $S$ we may assume that
$\hom(S',X)$ and $\hom(\L[m,j],X)$ are known to be manifolds.
Moreover, by $Kan(m,j)$, the bottom arrow is a submersion.
Therefore by transversality, $\hom(S,X)$ is a manifold.
\end{proof}

Since the horns $\L[m,j]$ are collapsible, we have:

\begin{cor}
Let $X$ be a Lie homotopy type, then $\hom(\L[m,j],X)$ is a
manifold.
\end{cor}

\subsection{The description of Lie 2-group(oid)s with finite data}

Often the conventional way with only finite layers of data to
understand Lie group(oid)s is more conceptual in differential
geometry. Therefore it is worth giving also the finite description
of Lie 2-groupoids. It turns out to be quite involved, however it
provides a direct generalization of Lie groupoids and it is useful
to pass to SLie groupoids. Moreover the version of 2-groupoid (not
Lie) is completely analogous to it.

The equivalent description of Lie 2-groupoids uses only three
layers $X_2 \Rrightarrow X_1 \Rightarrow X_0$ and associative
``3-multiplications'' \cite{duskin}. Following the notion of
simplicial manifolds, we call $d^n_i$ and $s^n_j$ the face and
degeneracy maps between $X_i$'s; they still satisfy the coherence
condition in \eqref{eq:face-degen}. To simplify the notation and
match it with the definition of groupoids, we use the notation
$\bt$ for $d^1_0$, $\bs$ for $d^1_1$ and $e$ for $s^0_0$. Then we
can safely omit the upper indices for $d^2_i$'s and $s^1_i$'s.
Actually we will omit the upper indices whenever it does not cause
confusion. Similarly to the horn spaces $\hom(\L[m,j], X)$, given
only these three layers, we define $\Lambda(X)_{m, j}$ to be the
space of $m$ elements in $X_{m-1}$ gluing along elements in
$X_{m-2}$ to a horn shape without the $j$-th face. \vspace{.6cm}

\centerline{\epsfig{file=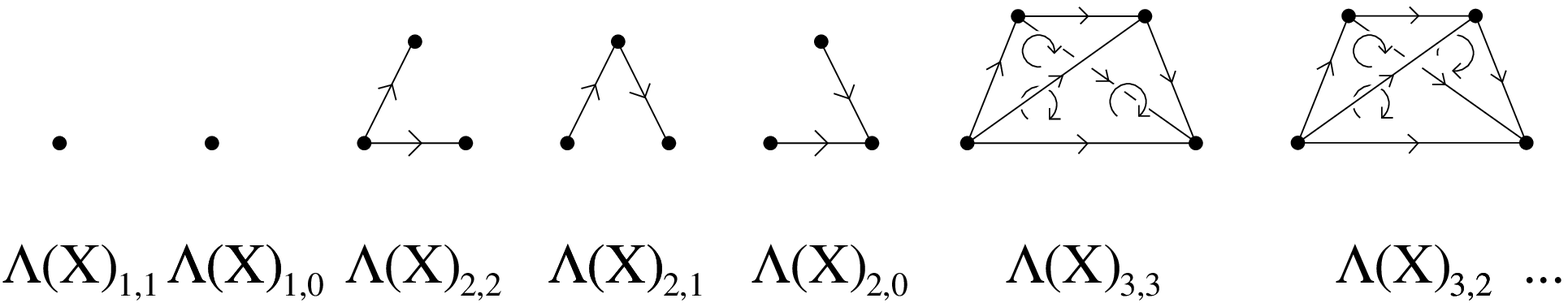,height=1.9cm}} \vspace{.6cm}
\noindent Here one imagines each $j$-dimensional face as an
element in $X_j$. For example,
\[
\begin{split}
\L(X)_{2,2}=X_1\times_{\bs, X_0, \bs}X_1, \; \Lambda(X)_{2, 1} &=
X_1 \times_{\bt, X_0, \bs} X_1, \;
\L(X)_{2,0}=X_1\times_{\bt, X_0, \bt} X_1, \\
\dots, \Lambda(X)_{3,0} &= (X_2\times_{d_2, X_1, d_1}
X_2)\times_{d_1\times d_2, \L(X)_{2,0}, d_1\times d_2} X_2.
\end{split}
\]

We remark that item \eqref{itm:st} and \eqref{itm:dd} in
the proposition-definition below imply that 
$\L(X)_{2,j}$'s and $\L(X)_{3,j}$'s  are manifolds. Then with this
condition we can define {\em 3-multiplications} as smooth maps
$m_i:\Lambda(X)_{3, i}   \to X_2$, $ i=0,\dots,3$. With
3-multiplications, there are natural maps between $\L(X)_{3,j}$'s.
For example, \[\L(X)_{3,0} \to \L(X)_{3,1}, \quad \text{ by}\;
(\eta_1, \eta_2, \eta_3) \to (m_0(\eta_1, \eta_2, \eta_3), \eta_2,
\eta_3). \] It is reasonable to ask them to be isomorphisms. In
fact, set theoretically, this simply says that the following four
equations are equivalent to each other:
\[
\begin{split}
\eta_0= m_0 (\eta_1, \eta_2, \eta_3), \quad \eta_1 =m_1(\eta_0,
\eta_2, \eta_3), \\ \eta_2=m_2(\eta_0, \eta_1, \eta_3), \quad
\eta_3=m_3(\eta_0, \eta_1, \eta_2).
\end{split}
\]

\begin{pdef}\label{def:finite-2gpd}
A Lie 2-groupoid can be also described by three layers $X_2
\Rrightarrow X_1 \rra X_0 $ and the following data:
\begin{enumerate}
\item  the face and degeneracy maps $d^n_i$ and $s^n_i$ satisfying
\eqref{eq:face-degen} for $n=1,2$ as explained above, and
\begin{enumerate}
\item\label{itm:st} [1-Kan] $\bt$ and $\bs$ are surjective submersions;
\item\label{itm:dd} [2-Kan] $d_0\times d_2: X_2\to \L(X)_{2,1}=X_1\times_{\bt, X_0, \bs} X_1$,
$d_0\times d_1: X_2 \to \L(X)_{2,2}=X_1\times_{\bs, X_0, \bs}X_1$,
and $d_1\times d_2: X_2\to \L(X)_{2,0}=X_1\times_{\bt, X_0, \bt}
X_1$ are surjective submersions.
\end{enumerate}
\item smooth maps (3-multiplications),
\[
m_i:\;  \Lambda(X)_{3, i}   \to X_2 \quad i=0,\dots,3.
\]such that
\begin{enumerate}
\item\label{itm:m-iso} the induced morphisms (by $m_j$ as above) $\L(X)_{3,i}\to \L(X)_{3,j}$ are all
isomorphisms;
\item $m_i$'s are compatible with the face and degeneracy maps:
\begin{equation}\label{coco}
\begin{split}
\eta=m_1(\eta, s_0\circ d_1(\eta), s_0\circ d_2(\eta)) & \big( \text{which is equivalent to} \; \eta=m_0(\eta, s_0 \circ d_1(\eta), s_0 \circ d_2 (\eta)) \big),\\
\eta= m_2(s_0 \circ d_0 (\eta), \eta, s_1 \circ d_2(\eta)) & \big( \text{which is equivalent to} \; \eta= m_1(s_0 \circ d_0 (\eta), \eta, s_1 \circ d_2(\eta)) \big),\\
\eta= m_3 (s_1 \circ d_0 (\eta), s_1 \circ d_1(\eta), \eta) &
\big( \text{which is equivalent to}\;  \eta= m_2 (s_1 \circ d_0
(\eta), s_1 \circ d_1(\eta), \eta) \big).
\end{split}
\end{equation}
\item $m_i$'s are associative, that is for a 5-simplex $(0, 1, 2, 3, 4)$,
\begin{equation}\label{pic:5-gon}
\xymatrix{
& 0 \ar @{<-}[ddl] \ar @{<-}[dddrr] \ar @{<-}[dr] \ar @{<-}[d] & &  \\
& 4\ar@{->}[dl]\ar@{->}[ddrr]& 1 \ar @{<-}[dll] \ar @{<-}[ddr] \ar @{<-}[l] & \\
3 \ar@{->}[drrr]  & & & \\
& & & 2 }
\end{equation}
if we are given faces $(0, i, 4)$'s and $(0, i, j)$'s in $X_2$,
where $i, j \in \{ 1, 2, 3\}$, then the following two methods to
determine the face $(1, 2, 3)$ give the same element in $X_2$:
\begin{enumerate}
\item $(1, 2, 3)= m_0( (0, 2, 3), (0,1, 3), (0, 1, 2))$;
\item we first obtain $(i, j, 4)$ using $m_i$'s on $(0, i, 4)$'s, then we have
\[(1, 2, 3)= m_3( (2, 3, 4), (1, 3, 4), (1, 2, 4)).\]
\end{enumerate}
\end{enumerate}
\end{enumerate}
\end{pdef}
\begin{remark}
Here we notice that the surjectivity of the maps in the 2-Kan
condition \eqref{itm:dd} insure the existence of  a usual (2-)
multiplication $m: X_1\times_{\bt, X_0, \bs} X_1 $ and an inverse
$i: X_1 \to X_1$ as explained in the introduction. For example
since $d_0\times d_2: X_2\to \L(X)_{2,1}$ is surjective  we could
take a section of it  and compose it with $d_1: X_2 \to X_1$. This
gives a multiplication map $m$, which is furthermore smooth if we
take a smooth section. This can always be achieved by passing to
an equivalent Lie 2-groupoid (see Section \ref{sec:equi-lht}).
However it is usually not associative on the nose, but only up to
elements in $X_2$. One could make it strictly associative by
passing to an infinite dimensional model.  But for us the
important example coming from integrating Lie algebroids does not
satisfy associativity on the nose. As stated in Theorem
\ref{thm:2-a}, the Lie 2-groupoid has $X_1=Lmor(T\Delta^1, A)$,
which is the space of $A$-paths satisfying the boundary condition
$a(0)=a(1)=0$. The 2-multiplication $m$ is simply the
concatenation of $A$-paths which is only associative up to
$A$-homotopies. Moreover we can not simply mod out $A$-homotopies
from $X_1$ since the resulting space might not be a smooth
manifold anymore.

On the other hand, only having a usual 2-multiplication $m$ and an
inverse map $i$, it is not guaranteed that the maps in the 2-Kan
condition \eqref{itm:dd} are submersions even when $m$ and $i$ are
smooth. But being submersions is in turn very important to prove
that $X_{n\geq 3}$ are smooth manifolds. Hence in the
differentiable category, we can not replace the 2-Kan condition by
a usual 2-multiplication and an inverse.
\end{remark}

\subsubsection*{The nerve of $X_2
\Rrightarrow X_1 \Rightarrow X_0$} To show that what we defined
just now is the same as Def. \ref{defngroupoids}, we form the {\em
nerve} of a Lie 2-groupoid $X_2 \Rrightarrow X_1 \Rightarrow X_0$
in Prop-Def. \ref{def:finite-2gpd}. We first define
\[X_3=\{ (\eta_0, \eta_1, \eta_2, \eta_3): \eta_0=m_0(\eta_1, \eta_2, \eta_3),
(\eta_1, \eta_2, \eta_3)\in \Lambda(X)_{3,0}\}. \] Then $X_3 \cong
\Lambda(X)_{3,0}$ is a manifold. Moreover, we have the obvious
face and degeneracy maps between $X_3$ and $X_2$,
\[
\begin{split}
& d^3_i( \eta_0, \eta_1, \eta_2, \eta_3)=\eta_i, \; i=0,\dots,3 \\
& s^2_0(\eta)=(\eta, \eta, s_0 \circ d_1 (\eta), s_0 \circ d_2 (\eta) ), \\
& s^2_1(\eta)=(s_0\circ d_0(\eta), \eta, \eta, s_1\circ d_2(\eta)), \\
& s^2_2(\eta)=(s_1\circ d_0(\eta), s_1\circ d_1(\eta), \eta,
\eta).
\end{split}
\]
The coherency \eqref{coco} insures that $s^2_i(\eta)\in X_3$. It
is also not hard to see that these maps together with $d^{\leq
2}_i$'s and $s^{\leq 1}_i$'s satisfy \eqref{eq:face-degen} for
$n\leq 3$.

Then the nerve can be easily described as the simplicial manifold
\[ X=cosk_2 (sk_2 (X_3\to X_2 \to X_1 \to X_0)). \]
More concretely, $X_n$ is made up by those $n$-simplices whose
2-faces are elements of $X_2$ and such that each set of four
2-faces gluing together as a 3-simplex is an element of $X_3$.
That is
\[ X_n = \{ f\in \hom_2( sk_2(\Delta_n), X_2)| f\circ (d_0\times d_1 \times d_2 \times d_3)(sk_3(\Delta_n))\subset X_3\} ,\]
where $\hom_2$ denotes the homomorphisms restricted to the 0,1,2
level and $X_2$ is understood as the tower $X_2 \Rrightarrow
X_1\Rightarrow X_0$ with all degeneracy and face maps. Then there
are obvious face and degeneracy maps which naturally satisfy
\eqref{eq:face-degen}.

However what is nontrivial and special in the differentiable
category is that the associativity of $m_i$'s assures that $X_n$
is a manifold. We prove this by an inductive argument. Let
$S_j[n]=sk_2(\Lambda[n,j])$. It is the contractible simplicial set
whose sub-faces all contain the vertex $j$ and whose only
non-degenerate faces are of  dimension 0, 1 and 2. Then similarly
to Lemma \ref{collapsable1}, we show that $\hom_2(S_j[n], X_2) $
is a manifold. Since $S_j[n]$ is constructed by adding 0,1,2
dimensional faces, it is formed by the procedure
\[
\xymatrix{ S' \ar[r] & S \\
           \L[n , j] \ar[u] \ar[r] & \Delta[n]\ar[u]}
\]
with $n\leq 2$. The dual pull-back diagram shows that
$\hom_2(S_j[n], X_2)$ is a manifold by induction
\[
\xymatrix{ \hom_2(S', X_2)\ar[d] &  \hom_2(S, X_2) \ar[l] \ar[d] \\
           \hom_2(\L[n,j] , X_2) & \hom_2(\D[n], X_2) ,  \ar[l] }
\]
since $\hom_2(\D[n], X_2) \to \hom_2(\L[n,j], X_2) $ are
surjective submersions by item \eqref{itm:st} and \eqref{itm:dd}
in the last Prop-Def.

Next we use induction to show that $X_n = \hom_2( S_0[n], X_2)$.
Similarly we will have $X_n = \hom_2( S_j[n], X_2)$.  It is clear
that $\tilde{f}\in X_n$ restricts to $\tilde{f}|_{S_0[n]} \in
\hom_2(S_0[n], X_2)$. We only have to show that $f\in
\hom_2(S_0[n], X_2)$ extends uniquely to $\tilde{f}\in X_n$. It is
certainly true for $n=0, 1, 2, 3$ just by definition. Suppose
$X_{n-1}= \hom_2( S_0[n-1], X_2)$. Then to get $f\in
\hom_2(S_0[n], X_2)$ from $f'\in \hom_2( S_0[n-1], X_2)$, we add a
new point $n$ and $n-1$ new faces $(0, i, n)$, $i\in \{1, 2,
\dots, (n-1)\}$ and dye them red\footnote{More precisely, they are
the image of these under the map $f$.}. Using 3-multiplication
$m_0$, we can determine face $(i, j, n)$ by $(0, i, n)$, $(0, j,
n)$ and $(0, i,j)$ and dye these newly decided faces blue.  Now we
want to see that each four faces attached together are in $X_3$,
then $f$ is extended to $\tilde{f} \in X_n$. We consider various
cases:
\begin{enumerate}
\item if none of the four faces contains the vertex $n$, then by
the induction condition, they are in $X_3$.
\item if one of the four faces contains $n$, then
there are three faces containing  $n$, we again have two
sub-cases:
\begin{enumerate}
\item if those three faces contain only one blue face of the form
$(i, j, n)$, $i, j \in \{ 1, \dots, (n-1)\}$, then the four faces
must contain three red faces and one blue face. According to our
construction, these four faces are in $X_3$;
\item if those three faces contains more than one blue face, then they must contain exactly three blue faces.  Then according to
associativity (inside the 5-gon $(0, i, j, k, n)$), these four
faces are also in $X_3$.
\end{enumerate}
\end{enumerate}
Now we finish the induction, hence $X_n$ is a manifold and it is
determined by the first three layers. So we have $\hom(\L[n,j],
X)=\hom_2(sk_2(\L[n,j]), X_2) =X_n$ and,

\begin{prop}The nerve $X$ of a Lie 2-groupoid $X_2\Rrightarrow X_1 \Rightarrow X_0$ as in
Prop-Def. \ref{def:finite-2gpd} is a Lie 2-groupoid as in Def.
\ref{defngroupoids}.
\end{prop}

\begin{prop}
The first three layers of a Lie 2-groupoid as in Definition
\ref{defngroupoids} is a Lie 2-groupoid as in Prop-Def.
\ref{def:finite-2gpd}.
\end{prop}
\begin{proof} The proof is more complicated and similar to the
case of 1-groupoids in the introduction. Here we point out that
the 3-multiplications $m_j$ are given by $Kan(3, j)$ and the
associativity is given by $Kan!(3, 0)$ and $Kan(4, 0)$.
\end{proof}

\section{SLie groupoids}
First we give the precise definition of SLie groupoids and
W-groupoids. This generalizes and completes the notion of
Weinstein (W-) groupoids in \cite{tz}. For example, we add some
new axioms on the level of 2-morphisms, and on the other hand,
find that some other axioms could be replaced or simplified.  Then
we point out some direct implications from the definition. The
notion of stacks has been extensively studied in algebraic
geometry for the past few decades. However stacks can also be
defined over other categories, such as the category of topological
spaces and category of smooth manifolds (see for example
\cite{SGA4} \cite{pronk} \cite{v2} \cite{bx1} \cite{metzler}). We
refer the readers to the latter two references for the concepts we
use here, such as (\'etale) differentiable stacks, their fibre
products, immersion maps between them. But for {\em surjective
submersion}, we adopt the definition in \cite{tz} which is a bit
different. $f: \cX \to \cY$ is a submersion if $X\times_{\cY} Y
\to Y$ is a submersion where $X$ and $Y$ are charts of $\cX$ and
$\cY$ respectively. $f$ is further a surjective submersion if it
is an epimorphism of differentiable stacks.

\begin{defi} [SLie (respectively W-) groupoid]\label{def:sliegpd} A  stacky Lie  or SLie (respectively Weinstein or W-) groupoid
over a manifold $M$ consists of the following data:
\begin{enumerate}
\item a differentiable (respectively an \'etale differentiable) stack $\cG$;
\item (source and target) maps $\bar{\bs}$,
$\bar{\bt}$: $\cG \to M$ which are surjective submersions between
differentiable stacks;
\item (multiplication) a map $m$: the fibre product $\cG\times_{\bbs, M, \bbt} \cG \to
\cG$, satisfying the following properties:
\begin{enumerate}
  \item\label{itm:m1} $\bbt \circ m=\bbt\circ pr_1$, $\bbs \circ m=\bbs\circ
  pr_2$, where $pr_i: \cG \times_{\bbs,M, \bbt} \cG \to \cG$ is the
  $i$-th projection $\cG\times_{\bbs, M, \bbt} \cG \to
\cG$;
  \item\label{itm:m-a} associativity up to a 2-morphism, i.e. there is a 2-morphism $a$ between
maps $m\circ (m \times id)$ and $m\circ(id\times m)$;
  \item\label{itm:a-higher} the 2-morphism $a$ satisfies a higher coherence described as following:
let the 2-morphisms on the each face of the
     cubes be $a_i$\footnote{All the $a_i$'s are generated by $a$, except that $a_4$ is $id$.} arranged in the following way:

front face (the one with the most $\cG$'s) $a_1$,  back $a_5$; up
$a_4$, down $a_2$; left $a_6$, right $a_3$,
$$
\xymatrix@=5pt{
     & & \cG\mathop\times\limits_{M}\cG\mathop\times\limits_{M}\cG \ar[dr]^{m\times id} \ar[ddd]^{id\times m}& \\
  \cG\mathop\times\limits_{M}\cG\mathop\times\limits_{M}\cG\mathop\times\limits_{M}\cG \ar[urr]^{id\times id\times m} \ar[dr]^{m\times id\times id} \ar[ddd]_{id\times m\times id} & & & \cG\mathop\times\limits_{M}\cG\ar[ddd]^{m} \\
     & \cG\mathop\times\limits_{M}\cG\mathop\times\limits_{M}\cG\ar[urr]^{id\times m} \ar[ddd]^{m\times id} & & & \\
     & & \cG\mathop\times\limits_{M}\cG \ar[dr]^{m} & \\
  \cG\mathop\times\limits_{M}\cG\mathop\times\limits_{M}\cG \ar[urr]^{id\times m} \ar[dr]^{m\times id} & & & \cG \\
     & \cG\mathop\times\limits_{M}\cG \ar[urr]^{m} & & }
$$
We require \[ (a_6\times id )\circ(id\times a_2)\circ(a_1\times
id)=  (id\times a_5) \circ (a_4 \times id) \circ (id \times a_3).
\]
\end{enumerate}

\item  (identity section) a morphism (respectively an immersion)  $\bar{e}$: $M\to \cG$
such that
\begin{enumerate}
\item\label{itm:e-b}  the following identities
\[
m\circ ((\bar{e}\circ \bbt)\times id)=id, \,\,m\circ (id\times
(\bar{e}\circ\bbs) )=id,\] hold\footnote{In particular, by
combining with the surjectivity of $\bbs$ and $\bbt$, one has
$\bbs \circ \bar{e}= id$, $\bbt \circ \bar{e}= id$ on $M$. In fact
if $x=\bbt(g)$, then $\bar{e}(x)\cdot g = g$ and $\bbt \circ m =
\bbt \circ pr_1$ imply that $\bbt(\bar{e}(x))= \bbt (g) = x$. } up
to 2-morphisms $b_l$ and $b_r$. Or equivalently there are two
2-morphisms
\begin{alignat*}{2}
  m \circ (id \times \bar{e}) & \overset{b_r}{\to} pr_1 : \cG \times_{\bbs, M} M \to \cG,&
 \quad
 m \circ (\bar{e} \times id) & \overset{b_l}{\to} pr_2 : M \times_{M, \bbt} \cG \to \cG , \\
 g\bar{e} (y) & \to g & \quad    \bar{e}(x) g & \to g
\end{alignat*}
where $y=\bbs(g)$ and $x=\bbt(g)$.
\item\label{itm:br}
We require the composed 2-morphism below, with $y=\bbs(g_2)$,
\[
g_1 g_2   \xrightarrow{b^{-1}_r}   (g_1 g_2) \bar{e}(y)
\xrightarrow{a}  g_1(g_2 \bar{e}(y)) \xrightarrow{b_r}  g_1g_2
\] to be the identity.\footnote{Forming this in the language of differentiable geometry, we notice that $pr_1 \circ (m \times id) $ and $m \circ
(pr_1\times pr_2)$ are the same map from $\cG\times_M \cG\times_M
M $ to $ \cG$, but as the diagram indicates,
\begin{equation}\label{diag:br}
\xymatrix{{\cG\times_M \cG\times_M M} \ar[d]^{ pr_1\times pr_2}
\ar@<-1ex>[d]_{id \times (m\circ (id \times \bar{e}))}
\ar[rr]^{ m \times id} & & {\cG \times_M M} \ar[d]^{ m\circ(id\times \bar{e})} \ar@<-1ex>[d]_{ pr_1}\\
{\cG \times_M \cG} \ar[rr]^m & & {\cG}, }
\end{equation}
they are related also via a sequence of 2-morphisms:
\begin{equation}\label{eq:br}
 pr_1 \circ (m \times id) \xrightarrow{b_r^{-1}\odot id}
m\circ(id\times \bar{e}) \circ (m \times id) \xrightarrow{a}  m
\circ (id \times (m\circ (id \times \bar{e})) ) \xrightarrow{id
\odot (id \times b_r)}  m \circ (pr_1\times pr_2)
\end{equation}
We require  the composed 2-morphisms be $id$, that is
\[ (id \odot ( id \times  b_r)) \circ a \circ (b_r^{-1} \odot id) =id
.\]}
\item\label{itm:bl}similarly with $x=\bbt(g_1)$,
\[ g_1 g_2   \xrightarrow{b^{-1}_l}   \bar{e}(x)(g_1 g_2)
 \xrightarrow{a^{-1}} (\bar{e}(x)g_1)g_2  \xrightarrow{b_l}  g_1g_2
 \] is the identity;
\item\label{itm:bl-br} with $y =\bbs(g)$ and $x=\bbt(g)$,
\[
g\xrightarrow{b_l^{-1}}  \bar{e}(x) g \xrightarrow{b_r^{-1}}
(\bar{e}(x)g) \bar{e}(y) \xrightarrow{b_l}  g\bar{e}(y)
\xrightarrow{b_r} g, \]is the identity.
\end{enumerate}

\item (inverse) an isomorphism of differentiable stacks
$\bar{i}$: $\cG \to \cG$ such that, up to 2-morphisms, the
following identities
\[ m\circ (\bar{i}\times id \circ \Delta)=\bar{e}\circ\bbs, \;\;
m\circ (id\times\bar{i}\circ \Delta)=\bar{e}\circ \bbt,\] hold,
where $\Delta$ is the diagonal map: $\cG\to \cG\times\cG$.
\end{enumerate}
\end{defi}
\begin{remark}
A W-groupoid is simply an \'etale SLie groupoid. This definition
of W-groupoid is different from the one in \cite{tz} in two
aspects: one is that here we add various higher coherences on
2-morphisms which make the definition more restricting but still
allow the W-groupoids $\cG(A)$ and $\cH(A)$, which are the
integration objects of the Lie algebroid $A$ constructed in
\cite{tz}; for the other see Remark \ref{rmk:wgpd-unnecessary}. On
the other hand, we do not add higher coherences for the
2-morphisms involving the inverse map. This is because the inverse
map can be removed from the definition. See Section
\ref{sec:inverse}.

Moreover, we notice that the 2-morphisms $\xymatrix@1{ \bar{e}(x)
\cdot \bar{e} (x) \ar[r]^{b_l}_{b_r}& \bar{e} (x) }$ are the same
because they are basically 2-morphisms between morphisms on a
manifold $M$. With some patience, we can check that the list of
coherences on 2-morphisms given here generates all the possible
coherences on these 2-morphisms. We also notice that the cube
condition \eqref{itm:a-higher} is  the differential version of the
pentagon condition \[ \left[ ((gh)k)l \to (g(hk))l\to g((hk)l) \to
g (h(kl)) \right]= \left[ ((gh)k)l\to (gh)(kl)\to g(h(kl)) \right]
.\]
\end{remark}

\subsection{Good charts} \label{sec:embedding}
Given an SLie groupoid $\cG \rra M$, the identity map $\bar{e}:
M\to \cG$ corresponds to an Hilsum-Skandalis (H.S.) morphism from
$M\rra M$ to $G_1\rra G_0$ for some presentation of $\cG$. But it
is not clear whether $M$ embeds into $G_0$. It is not even obvious
whether there is a map $M\to G_0$. In general, one could ask: if
there is a map from a manifold $M$ to a differentiable stack
$\cG$, when can one find a chart $G_0$ of $\cG$ such that $M\to
\cG$ lifts to $M\to G_0$, namely when is the H.S. morphism $M\rra
M$ to $G_1\rra G_0$ a strict groupoid morphism? If the stack $\cG$
is \'etale, can we find an \'etale chart $G_0$? We answer these
questions by the following lemmas.  It turns out that the \'etale
case is easier and when $M\to \cG$ is an immersion we can always
achieve an \'etale chart.

\begin{lemma}\label{lemma:embedding}
For an immersion $\bar{e}: M\to \cG$ from a manifold $M$ to an
\'etale stack $\cG$,  there is an \'etale chart $G_0$ of $\cG$
such that $\bar{e}$ lifts to an embedding $e: M\to G_0$. We call
such charts {\em good} charts and their corresponding groupoid
presentations {\em good} presentations.
\end{lemma}
\begin{proof}
Take an arbitrary \'etale chart $G_0$ of $\cG$. The idea is to
find an ``open neighborhood'' $U$ of $M$ in $\cG$ with the
property that $M$ embeds in $U$ and there is an \'etale
representable map $U\to \cG$. Since $G_0\to \cG$ is an \'etale
chart, in particular epimorphic, $G_0\sqcup U \to \cG$ is an
\'etale representable epimorphism\footnote{Note that being \'etale
implies being submersive.}, that is, a new \'etale chart of $\cG$.
Then the lemma is proven since $M\hookrightarrow G_0\sqcup U$ is
an embedding.

Now we look for such a $U$. Since $\bar{e}:M \to \cG$ is an
immersion, the pull-back $M\times_{\cG}G_0 \to G_0$ is an
immersion and $M\times_{\cG}G_0 \to M$ is an \'etale epimorphism.
We cover $M$ by small enough open charts $V_i$'s so that $V_i$
lifts to isomorphic open charts $V'_i$ on $M\times_{\cG} G_0$.
Then $V'_i \to G_0$ is an immersion so locally it is an embedding.
Therefore we can divide $V_i$ into even smaller open charts
$V_{i_j}$ such that $V_{i_j}\cong V'_{i_j}\to G_0$ is an
embedding. Hence we might assume that the $V_i$'s form an open
covering of $M$ such that $\bar{e}$ lifts to embeddings $e_i: V_i
\hookrightarrow G_0$. It appears with the language of the
Hilsum-Skandalis (H.S.) bibundles as the diagram on the right,
\[
\xymatrix{V'_i\subset & M\times_{\cG} G_0 \ar[r] \ar[d] & G_0 \ar[d] \\
V_i \subset & M \ar[r] & \cG} \quad
\xymatrix{ M\ar[dd]\ar@<-1ex>[dd]&  &G_1\ar[dd]\ar@<-1ex>[dd]\\
& V'_i\subset M\times_{\cG}G_0\ar[dr]^{J_r}\ar[dl]_{J_l}& &\\
M\supset V_i \ar@<-1ex>[ur]_{\sigma_i}&  &G_0\\
}
\]
Here $e_i=\sigma_i\circ J_r$. Since the action of $G_1$ on the
H.S. bibundle is free and transitive, there exists a unique
groupoid bisection $g_{ij}$ such that $e_i \cdot g_{ij}=e_j$ on
the overlap $V_i\cap V_j$ (in fact $G_1|_{\sqcup
V_i}=\bt^{-1}(\sqcup e_i(V_i))\cap \bs^{-1}(\sqcup e_i(V_i)) \rra
\sqcup e_i(V_i)$ is Morita equivalent to $M$, where $\bt, \bs:
G_1\rra G_0$). Since $G_1 \rra G_0$ is \'etale, bisection $g_{ij}$
extends uniquely to $\bar{g}_{ij}$ on an open set
$\bar{U}_{ij}\subset G_1$. Moreover, there exist open sets $U_i
\supset e_i(V_i)$ of $G_0$ such that
\[e_i(V_i\cap V_j)\subset \bt(\bar{U}_{ij}) =: U_{ij} \subset U_i,
\quad e_j(V_i\cap V_j)\subset \bs(\bar{U}_{ij}) =: U_{ji} \subset
U_j. \] Since $e_j \cdot g_{ij}^{-1}=e_i$ these sets are well defined.   

\centerline{
\epsfig{file=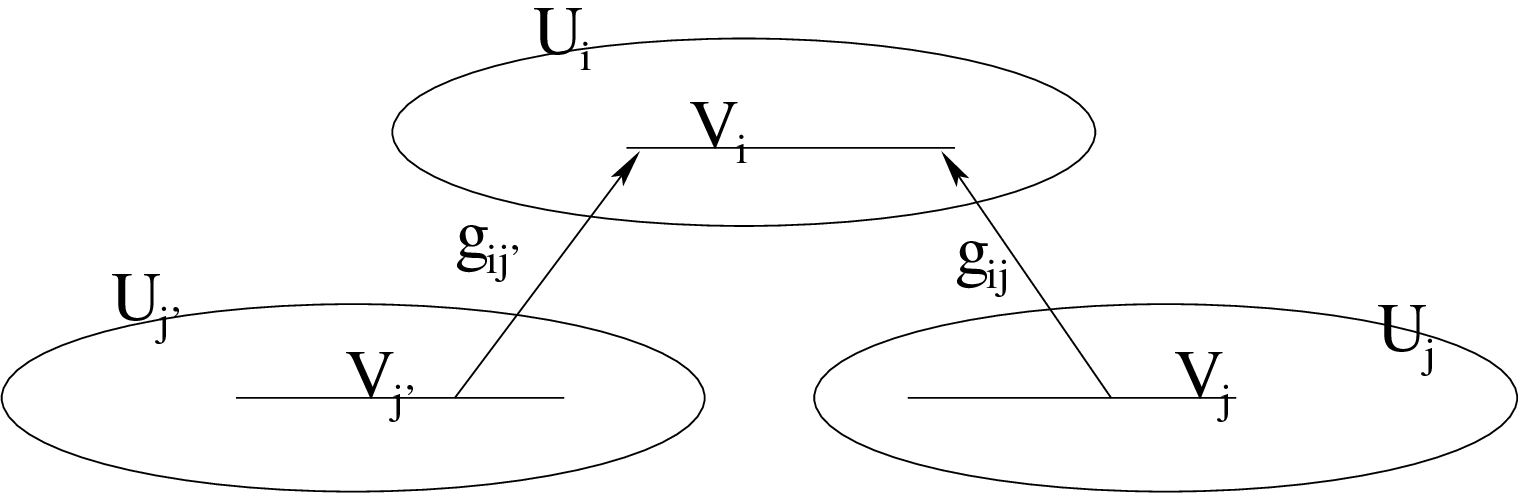,height=2cm}}

Because of uniqueness and $g_{ij}\cdot g_{jk}=g_{ik}$, we have
$\bar{g}_{ij}\cdot \bar{g}_{jk}=\bar{g}_{ik}$ on the open subsets
$\bar{U}_{ijk}:= \{ (\bar{g}_{ij}, \bar{g}_{jk}, \bar{g}_{ik}):$
when $\bar{g}_{ij}\cdot \bar{g}_{jk}$ exists and in
$\bar{U}_{ik}$.$\}$. Then \[e_i(V_i \cap V_j \cap V_k) \subset
U_{ijk}:= \bt (Im(\bar{U}_{ijk} \to \bar{U}_{ij})) \subset
U_{ij}\cap U_{ik} \subset U_i, \] and similarly for $j$ and $k$.
Therefore with these $U$'s we are in the situation of a {\em germ
of manifolds}  of $M$ defined as below.

A {\em germ of manifolds at a point} $m$ is a series of manifolds
$U_i$'s containing $m$ such that $U_i$ agrees with $U_j$ in a
smaller open set $(m\in) U_{ij}\subset U_i$ by $x\sim f_{ij}(x)$.
A {\em compatible riemannian metric} of a germ of manifolds
consists of a riemannian metric $g^i$ on each $U_i$ such that two
such riemannian metrics $g^i$ and $g^j$ on $U_i$ and $U_j$ agree
with each other in the sense that $g^i(x)=g^j(f_{ij}(x))$ in a
smaller open set (possibly a subset of $U_{ij}$).  With this one
can define the exponential map $exp$ at $m$ using the usual
exponential map of a riemannian manifold provided the germ is
finite, namely there are finite many manifolds in the germ (which
is true in our case, since $V_i$ intersects finite other $V_j$'s).
Then $exp$ gives a Hausdorff manifold containing $m$.

If a series of locally finite manifolds $U_i$'s and morphisms
$f_{ij}$'s form a germ of manifolds for every point of a manifold
$M$, we call it {\em a germ of manifolds of } $M$. Here local
finiteness means that any open set in $M$ is contained in finitely
many $U_i$'s and $M$ has the topology induced by the $U_i$'s, that
is $M\cap U_i$ is open in $M$.  We can always endow each of them
with a compatible riemannian metric, beginning with any riemannian
metric $g^i$ on $U_i$ and modifying it to the sum
$g'^i(x):=\sum_{k, x\in U_{ik}} g^k(f_{ik}(x))$ (with
$f_{ii}(x)=x$) at each point $x \in U_i$. In this situation, one
can take a tubular neighborhood $U$ of $M$ by the $exp$ map of the
germ. Then $U$ is a Hausdorff manifold.

Applying the above construction to our situation, we have a
Hausdorff manifold $U \supset M$ with the same dimension as $G_0$.
$U$ is basically glued by small enough open subsets
$\tilde{U}_i=U\cap
U_i$ containing $V_i$'s along $\tilde{U}_{ij}:=  U\cap U_{ij}$ 
so that the gluing result $U$ is still a Hausdorff manifold.
Therefore $U$ is presented by $\sqcup \tU_{ij} \rra \sqcup \tU_i$
which maps to $G_1\rra G_0$ via $U_{ij}\cong
\bar{U}_{ij}\hookrightarrow G_1$. So there is a map $\pi: U\to
\cG$. Since $\tU_{i} \to G_0$ are \'etale maps, by the technical
lemma below, $\pi$ is a representable \'etale map.
\end{proof}

\begin{lemma}\label{lemma:rep-sub}
Given a manifold $X$ and an (\'etale) differentiable stack $\cY$,
a map $f: X\to \cY$ is an (\'etale) representable submersion if
and only if there exists an (\'etale) chart $Y_0$ of $\cY$ such
that the induced local maps $X_i\to Y_0$ are (\'etale)
submersions, where $\{X_i\}$ is an open covering of $X$.
\end{lemma}
\begin{proof}
For any $V\to \cY$, $X_i\times_{\cY} V = X_i \times_{Y_0} Y_0
\times_{\cY} V $ is representable and $X_i\times_{\cY} V \to V$ is
an (\'etale) submersion since $X_i\to Y_0$ and $Y_0\to \cY$ are
representable (\'etale) submersions. Since $X_i$ glue together to
$X$, $X_i\times_{\cY} V$ with the inherited gluing maps glue to a
manifold $X\times_{\cY} V$. Since being an (\'etale) submersion is
a local property, $X\times_{\cY} V \to V$ is an (\'etale)
submersion.
\[
\xymatrix{ X_i \times_{\cY} V \ar[drr] \ar[dd] \\
& \curvearrowright & X_k \times_{\cY} V \ar[dll] \\
X_j \times_{\cY} V } \quad \xymatrix{&\\ \Leftarrow \\ &} \quad \xymatrix{  X_i  \ar[drr] \ar[dd] \\
& \curvearrowright & X_k  \ar[dll] \\
X_j  }
\]
\end{proof}

\begin{remark}\label{rk:loc} If $\bar{e}$ is the identity map of W-groupoid $\cG \rra M$, then an
open neighborhood of $M$ in $U$ has an induced local groupoid
structure from the stacky groupoid structure \cite[Section~5]{tz}.
\end{remark}
\begin{remark} \label{rmk:wgpd-unnecessary}
The following assumption in the original definition of W-groupoids
becomes obviously unnecessary after Lemma \ref{lemma:embedding}:
\\
\noindent {\scriptsize ``Moreover, restricting to the identity
section, the above 2-morphisms between maps are the $id$
2-morphisms. Namely, for example, the 2-morphism $\alpha$ induces
the $id$ 2-morphism between the following two maps:\[ m\circ ((m
\circ (\bar{e}\times\bar{e}\circ \delta))\times \bar{e} \circ
\delta)=m\circ(\bar{e}\times(m\circ(\bar{e}\times\bar{e}\circ\delta))\circ\delta),
\]where $\delta$ is the diagonal map: $M\to M\times M$.''}
\\
This is unnecessary because $(M\rra M)$ embeds into $(G_1\rra
G_0)$ on both levels of the groupoid  via $\bar{e}$ and the
restricted maps become maps between manifolds and therefore there
are no nontrivial 2-morphisms between them.
\end{remark}

We further prove the same lemma in the non-\'etale case.
\begin{lemma}\label{lemma:embedding-non-etale}
For a morphism $\bar{e}: M\to \cG$ from a manifold $M$ to a
differentiable stack $\cG$,  there is a chart $G_0$ of $\cG$ such
that $\bar{e}$ lifts to an embedding $e: M\to G_0$. We  call also
such groupoid presentations {\em good} presentations.
\end{lemma}
\begin{proof} We follow the proof of the \'etale case, but replace ``\'etale'' with ``submersion''.
We need a $U$ with a representable submersion  to $\cG$ and an
embedding of $M$ into $U$. There are two differences: first, $V_i$
embeds in $V'_i$ instead of being isomorphic to it, and we do not
have an embedding $V'_i \hookrightarrow G_0$; second, since
$G_1\rra G_0$ is not \'etale, the bisection $g_{ij}$ does not
extend uniquely to some $\bar{g}_{ij}$ and we can not have the
cocycle condition immediately.

The first difference is easy to compensate: given any morphism $f:
N_1 \to N_2 $, we can always view it as a composition of an
embedding and a submersion as $N_1 \overset{id\times
f}{\hookrightarrow} N_1\times N_2
\overset{pr_2}{\twoheadrightarrow} N_2$. In our case, we have the
following decomposition $M\times_{\cG}G_0 \hookrightarrow H_0
\twoheadrightarrow G_0$, then we use the pull-back groupoid $H_1:=
G_1\times_{G_0\times G_0} H_0\times H_0$ over $H_0$ to replace
$G$. Thus we obtain an embedding of $V'_i \to H_0$ and thus an
embedding $V_i \to H_0$. Then since $H_1\rra H_0$ is Morita
equivalent to $G_1\rra G_0$, we just have to replace $G$ by $H$ or
call $H$ our new $G$. It was not possible to do so in the \'etale
case since $H_0$ might not be an \'etale chart of $\cG$.

For the second difference, first of all we could assume $M$ to be
connected to construct such $U$. Otherwise we take the disjoint
union of $U$'s of each connected component of $M$.

Then take any $V_i$ and consider all the charts $V_j$'s
intersecting $V_i$. We choose $\bar{g}_{ij}$ extending $g_{ij}$ on
an open set $\bar{U}_{ij}$. As before we define the open sets
$U_i$, $U_j$'s, and $U_{ij}$. Then for $V_j$ and $V_{j'}$ both
intersecting $V_i$, we choose $\bar{g}_{jj'}$ to be the one
extending (see below) $\bar{g}_{ij}^{-1}\bar{g}_{ij'}$ with
$\bs(\bar{g}_{ij}^{-1}\bar{g}_{ij'})$ in the triple intersection
$\bg_{ij'}^{-1}\cdot(\bg_{ij}\cdot U_j ) \cap U_{j'}$ where
multiplication applies when it can. Since $\bar{g}$'s are local
bisections, $\bar{g} \cdot $ is an isomorphism. Identifying via
these isomorphisms, we view and denote the above intersection as
$U_{j'ij}$ for simplicity.

Now we clarify in which sense and why the extension always exists.
Let us assume $\dim M =m$, $\dim G_i =n_i$. 
Here we identify $V_j$ with its embedded image in $G_0$. 
Then since we are dealing with local charts, we might assume that
both $\bt$ and $\bs$ of $G_1\rra G_0$ are just projections from
$\R^{n_1}$ to $\R^{n_0}$. A section of $\bs$ is a vector valued
function $\R^{n_0}\to \R^{n_1}/\R^{n_0}$, and it being  a
bisection, namely also a section of $\bt$,  is an open condition.
We can always perturb a section to get a bisection. If we can
extend $\bg_{ij}^{-1}\bg_{ij'}$ and $g_{jj'}$ from $U_{j'ij} \cup
e_{j'}(V_j\cap V_{j'})$ to a bisection $\bg_{jj'}$ such that $\{
\bs(\bar{g}_{jj'})\}$ is an open set in $U_{j'}$, then we obtain a
bisection $\bar{g}_{jj'}$ from $U_{j'j}:=
\bg_{jj'}^{-1}(\{\bt(\bg_{jj'})\}\cap U_{j})$ to $U_{jj'} :=
\{\bt(\bg_{jj'})\}\cap U_{j}$. It is easy to see that
$U_{jj'}\cong U_{j'j}$ are open in $G_0$ since
$\{\bt(\bg_{jj'})\}\cong \{\bs(\bg_{jj'})\}$.

Therefore we are done as long as we can extend a smooth function
$f$ from the union of an open submanifold $O$ with a closed
submanifold $V$ of an open set $B\subset \R^{n_0}$ to the whole
$B$. Since $V$ is closed, using its tubular neighborhood and
partition of unity, we can first extend $f$ from $V$ to $B$ as
$\tilde{f}$. Then $f_1= f-\tilde{f}|_{O\cup V}$ is 0 on $V$. We
shrink the open set $O$ a little bit to $O_i$ such that $V \cap O
\subset O_2\subset O_1 \subset O$. Then we always have a smooth
function $p$ on $B$ with $p|_{\bar{O}_2 } =1$ and $p|_{B-O_1} =0$.
Then the extension function $\tilde{f}_1$ is defined by
\[\tilde{f}_1(x)=
  \begin{cases}
    f_1(x) \cdot p(x) & \text{$x\in O$}, \\
    0 & \text{otherwise}.
  \end{cases}
\] It is easy to see that $\tilde{f}_1$ is smooth and it agrees with $f_1$ on $O_2$ and $V$ because
$V-O_2=V-O_1 \subset B-O_1$ and $p|_{V-O_2}=0$. Hence $\tilde{f} +
\tilde{f}_1$ extends $f|_{O_2\cup V}$. Now we extend the
$\bar{g}_{ij}^{-1}\bar{g}_{ij'}$'s to $\bar{g}_{jj'}$'s, then the
$\bar{g}$'s satisfy the cocycle condition on smaller open sets of
the triple intersections $U_{j'ij}$ by construction.

Then we view $V_i \cup(\cup_{j:V_i\cap V_j \neq \emptyset} V_j)$
as one chart. Notice that a connected manifold is path connected.
Also notice that we didn't use any topological property of $V_i$
or $U_i$. This construction will eventually extend to the whole
manifold $M$ and obtain the desired $\bg_{ij}$'s. Therefore we are
again in the situation of a germ of manifolds and we can apply the
proof of Lemma \ref{lemma:embedding} to get the result.
\end{proof}

\section{From SLie-groupoids to Lie 2-groupoids}\label{sec:slie-2}

Suppose $\cG\rra M$ is an SLie groupoid, in this section we
construct a Lie 2-groupoid $X_2 \Rrightarrow X_1 \Rightarrow X_0$
corresponding to it. When $\cG\rra M$ is a W-groupoid, the
corresponding Lie 2-groupoid is {\em 2-\'etale}, that is the maps
$X_2 \to \hom(\Lambda[2,j], X)$ are \'etale for $j=0,1,2$.

\begin{thm}\label{thm:slie-2} An  SLie (respectively W-) groupoid $\cG \rra M$ with a chosen good chart
(respectively good \'etale chart) $G_0$ of $\cG$ corresponds to a
Lie 2-groupoid (respectively 2-\'etale Lie 2-groupoid) $X_2
\Rrightarrow X_1 \Rightarrow X_0$.
\end{thm}

\subsection{The construction of $X_2 \Rrightarrow X_1 \Rightarrow
X_0$}\label{sec:statement} Given an SLie groupoid $\cG\rra M$, let
$G_1 \rra G_0$ be a good groupoid presentation of $\cG$ and $E_m$
a bimodule presenting the morphism $m$. Let $J_l: E_m \to
G_0\times_M G_0$ and $J_r: E_m \to G_0$ be the moment maps of the
bimodule $E_m$. Notice that for an SLie groupoid $g\cdot 1 $``=''1
up to a 2-morphism, that is $m|_{\cG \times_M M}=id$ up to a
2-morphism. Translating this into groupoid language,  $J_1^{-1}
(G_0\times_M M)$ is the bimodule presenting $m|_{\cG \times_M M}$.
By the definition of SLie groupoids, there is an isomorphism $b_r:
J_l^{-1}(G_0\times_M M)\to G_1$. Similarly, there is an
isomorphism $b_l: J_l^{-1}(M \times_M G_0) \to G_1$.

We construct
\[ X_0 =M , X_1=G_0, X_2=E_m \]
with the structure maps
\begin{equation} \label{eq:stru-maps}
\begin{split}
 \d^1_0=\bs, \d^1_1=\bt: X_1\to X_0, \quad &\d^2_0=pr_2 \circ J_l, \d^2_1=J_r, d^2_2=pr_1\circ J_l  : X_2 \to X_1,\\
 s^0_0=e: X_0 \to X_1, \quad
 & s^1_0= b_l^{-1}\circ e_G, s^1_1= b_r^{-1}\circ e_G: X_1 \to X_2
\end{split}
\end{equation} where $pr_i$ is the i-th projection
$G_0\times_M G_0 \to G_0$, $\bs$, $\bt$ present $\bbs$, $\bbt$:
$\cG \rra M$, and $e_G$ is the identity embedding $G_0 \to G_1$.
We still need the 3-multiplication maps
\[
m_i:\;  \Lambda(X)_{3, i}   \to X_2 \quad i=0,\dots,3.
\]
Let us first construct $m_0$. Notice that in the 2-associative
diagram, we have a 2-morphism $a: m\circ(m \times id) \to m\circ
(id \times m)$. Translating this into the language of groupoid, we
have the following isomorphism of bimodules:
\begin{equation}\label{eq:a}
a: (( E_m \times_{G_0} G_1)\times_{G_0 \times_M G_0}E_m )/ (G_1
\times_M G_1) \to ( (G_1 \times_{G_0} E_m)\times_{G_0 \times_M
G_0}E_m)/(G_1 \times_M G_1).
\end{equation}
Suppose $(\eta_1, \eta_2, \eta_3) \in \Lambda(X)_{3,0}$. Then
$(\eta_3, 1, \eta_1)$ represents a class in $(E_m\times_{G_0} G_1)
\times_{G_0 \times_M G_0}E_m/ \sim$ (we write $\sim$ when it is
clear which groupoid action  is meant). Moreover, its image under
$a$ can be represented by $( 1, \eta_0,\eta_2)$, that is,
\[ a( [(\eta_3, 1,\eta_1)] )=  [( 1, \eta_0,\eta_2)].
\quad
\begin{xy}
*\xybox{(0,0);<3mm,0mm>:<0mm,3mm>::
  ,0
  ,{\xylattice{-5}{0}{-4}{0}}}="S"
  ,{(-10,-10)*{\bullet}}, {(-10, -12)*{_1}},
     ,{(0,0)*{\bullet}}, {(0, 2)*{^{0}}}, {(10, -10)*{\bullet}},
     {(10, -12)*{^{2}}}, {(15, -4)*{\bullet}}, {(17,-5)*{^3}},
     {(-10, -10) \ar@{->}^{g_1} (0,0)},
     { (10, -10) \ar@{->}^{g_2} (-10,-10)},
      { (10, -10) \ar@{->} (0,0)},
     {(15, -4)\ar@{->}^{g_3} (10, -10)},
     {(15, -4)\ar@{->}(0, 0)},
     {(15, -4)\ar@{->}(-10,-10)},
\end{xy}  \]
As before we imagine the $j$-dimensional faces of the picture are
elements of $X_j$. Then we arrive naturally at $\eta_0$. To show
that the above definition is good, we have to show: 1) the image
under $a$ can be represented by some element of the form $( 1,
\eta_0,\eta_2)$; 2) the choice of $\eta_0$ is unique; 3) this map
is smooth.

Suppose $a([ \eta_3, 1,\eta_1])= [( 1, \teta_0,\teta_2)]$. Notice
that the $G_1\times_M G_1$ action on $G_1$ on both sides of
\eqref{eq:a} are the right multiplication by one of the copies of
$G_1$, therefore we can always suppose the element in $G_1$ is 1.

Notice that similarly to the proof of Lemma \ref{proof-morita},
the map $m \times pr_1: \cG \times_{\bbs, M, \bbt} \cG \to \cG
\times_{\bbt, M, \bbt}\cG$ is an isomorphism. Therefore $E_m
\times_{pr_1 \circ J_l, G_0, \bt_G} G_1$ is a Morita bibundle from
the Lie groupoid $G_1 \times_{\bs\circ \bs_G, M, \bt\circ \bs_G}
G_1\rightrightarrows G_0 \times_{\bs , M , \bt} G_0$ to
$G_1\times_{\bt\circ \bs_G, M, \bt \circ \bs_G}G_1
\rightrightarrows G_0\times_{\bt, M , \bt} G_0$, where $\bs_G$,
$\bt_G$ are the source and target maps on $G$ and $\bs$, $\bt$ are
the maps $G_0\to M$ presenting $\bbs$, $\bbt$. Here the two moment
maps are $J_l$ (for $E_m$) and $J_r \times \bs_G$. Therefore the
left groupoid action of $G_1 \times_{\bs\circ \bs_G, M, \bt\circ
\bs_G} G_1$ is principal on the bibundle $E_m \times_{pr_1 \circ
J_l, G_0, \bt_G} G_1$. Notice that  $(\eta_2, 1)$ and $(\teta_2,
1)$ are on the same fibre of $J_r\times \bs_G$ over
$G_0\times_{\bt, M , \bt} G_0$ because the way we arrange the map
makes them share the same edges $1\to 0$ and $3\to 0$. So there is
a unique groupoid element $(\gamma_1, \gamma_2)$ such that
$(\teta_2, 1) \cdot (\gamma_1, \gamma_2) = (\eta_2, 1)$. This also
forces $\gamma_2$ to be 1. Therefore, we have
\[ (\teta_2, 1,\teta_0)\cdot (\gamma_1, 1)= (\eta_2,1, \eta_0). \]
By uniqueness of $\gamma_1$, the choice of $\eta_0$ is unique.
Moreover, since smoothness is a local property and the map $a$ is
an isomorphism and the groupoid action is principal, our map is
naturally smooth.

For other $m$'s, we can precede in a similar fashion.  More
precisely, for $m_1$ one can make the same definition for $m_0$
but using $a^{-1}$. It is even easier to define $m_2$ and $m_3$.
For example, for $m_2$, on the right hand side of \eqref{eq:a},
since the $G_1 \times_M G_1$ action on $G_1 \times_M E_m$ is
principal, the analogous statement to 1) naturally holds. Thus we
realize that given any three $\eta$'s, we can always put them in
the same spots as we did for $m_0$. Then any three of them
determine the fourth. Hence the $m$'s are compatible with each
other.

\subsection{Proof that what we construct is a Lie 2-groupoid}

By Prop-Def. \ref{def:finite-2gpd}, to show the above construction
gives us a Lie 2-groupoid, we just have to show that the $m$'s
satisfy the coherence conditions and associativity, and 1-Kan
2-Kan conditions. Condition 1-Kan is implied by the fact that
$\bs,\bt: G_0\rra M$ are surjective submersions; $Kan(2,1)$ is
implied by the fact that the moment map $J_l: E_m\to
G_0\times_{\bs,M, \bt, } G_0$ is a surjective submersion;
$Kan(2,0)$ is implied by the surjective submersion
\[J_l: E_m\cong \big( G_1\times_{M}
E_i\times_{G_0\times_M G_0} E_m \big)/G_1\times_M G_1 \to
G_0\times_{\bt, M, \bt} G_0 =\Lambda[2,0](X), \] where $E_i\cong
G_1$ (Lemma \ref{isom}) is the bimodule presenting the inverse map
of $\cG$, the composed bibundle is the one presenting the map
$(g,h)\mapsto (g,h^{-1}) \mapsto gh^{-1}$; $Kan(2,2)$ follows
similarly as $Kan(2,0)$.

\subsubsection*{The coherence conditions}
The first identity in \eqref{coco} corresponds to an identity of
2-morphisms,
\[ \big (1\cdot (g_1 \cdot g_2) \overset{a}{\sim} (1\cdot g_1)\cdot g_2
\sim g_1\cdot g_2 \big) = \big( 1\cdot (g_1 \cdot g_2) \sim
g_1\cdot g_2 \big),
\] Restrict the two bimodules in \eqref{eq:a} to $M \times_M G_0 \times_M
G_0$, then we get $E_m$ on the left hand side because $J_l^{-1} (M
\times_M G_0)=G_1$ and $\big( (G_1\times_M G_1)
\times_{G_0\times_M G_0}E_m \big) /G_1\times_M G_1 =E_m$. More
precisely, the elements in $(E_m
\times_{G_0}G_1)\times_{G_0\times_M G_0} E_m|_{M\times_M G_0
\times_M G_0} /\sim$ have the form $[(s_0\circ d_2 (\eta), 1,
\eta)]$, and  the isomorphism to $E_m$ is given by $[(s_0\circ d_2
(\eta), 1,\eta)]\mapsto \eta$. Similarly for the right hand side,
i.e. $[(s_0\circ d_1(\eta), 1, \eta)]\mapsto \eta$ gives the other
isomorphism. By \ref{itm:bl} in Def. \ref{def:sliegpd}, the
composition of the first and the inverse of the second map is $a$
(restricted on the restricted bimodules), so we have
\[ a ([(s_0\circ d_2 (\eta), 1,\eta )])= ([(1, \eta,s_0\circ
d_1(\eta))]),
\]
which implies the first identity in \eqref{coco}. The rest follows
similarly.

\subsubsection*{Associativity} \label{sec:3-asso}
For the associativity we will have to use the cube condition
\ref{itm:a-higher} for $a$ in Def. \ref{def:sliegpd} (called also
``pentagon condition'' in the literature). Let $\eta_{ijk}$ denote
the faces in $X_2$ fitting in diagram \eqref{pic:5-gon}. Suppose
we are given the faces $\eta_{0i4}'s\in X_2$ and  the faces
$\eta_{0ij}'s\in X_2$. Then we have two ways to determine the face
$\eta_{123}$ using $m$'s as described in Prop-Def.
\ref{def:finite-2gpd}. We will show below that these two
constructions give the same element in $X_2$.

Translate the cube  condition into the language of Lie groupoids.
The morphisms become bimodules and the 2-morphisms become the
morphisms between bimodules. The cube condition tells us that the
following two compositions of morphisms are the same (here for
simplicity, we omit writing the base space of the fibre products
and the groupoids by which we take quotients):
\[
\begin{split}
     &  (E_m \times G_1 \times G_1)\times (E_m \times G_1) \times E_m /\sim  \quad <--> \quad ((g_1g_2)g_3)g_4 \\
\overset{id\times a}{\lra}
& (E_m \times G_1 \times G_1)\times (G_1 \times E_m) \times E_m /\sim   \quad <--> \quad (g_1 g_2)(g_3g_4)\\
\overset{id}{\lra}
& (G_1\times G_1 \times E_m) \times (E_m \times G_1)\times E_m/\sim  \quad <--> \quad (g_1 g_2) (g_3 g_4) \\
\overset{id \times a} {\lra} & (G_1 \times G_1 \times E_m)\times
(G_1 \times E_m )\times E_m /\sim  \quad <--> \quad g_1(g_2(g_3
g_4))
\end{split}
\]
and
\[
\begin{split}
  & (E_m \times G_1 \times G_1)\times (E_m \times G_1) \times E_m /\sim  \quad <--> \quad ((g_1g_2)g_3)g_4 \\
\overset{a\times id}{\lra} &(G_1\times E_m \times G_1) \times (E_m
\times G_1) \times E_m /\sim  \quad <--> \quad   (g_1(g_2g_3))g_4\\
\overset{id\times a }{\lra}
& (G_1 \times E_m \times G_1) \times (G_1 \times E_m)\times E_m/\sim  \quad <--> \quad  g_1((g_2g_3)g_4)\\
\overset{ a \times id} {\lra} & (G_1 \times G_1 \times E_m) \times
(G_1 \times E_m )\times E_m/\sim  \quad <--> \quad
g_1(g_2(g_3g_4))
\end{split}
\]
Tracing through where the element $(\eta_{034}, (\eta_{023}, 1),
(\eta_{012}, 1,1))$ goes via the first and second composition, it
should end up in the same element. So we have
\[
\begin{split}
  & [((\eta_{012}, 1, 1)), (\eta_{023}, 1),\eta_{034} ] \\
\overset{id \times a}{\mapsto}
& [ ( (\eta_{012}, 1, 1), (1, \eta_{234}),\eta_{024} )] \\
\overset{id}{\mapsto}
& [ ((1,1, \eta_{234}), (\eta_{012}, 1), \eta_{024}) ] \\
\overset{id \times a}{\mapsto} & [((1,1, \eta_{234}), (1,
\eta_{124}), \eta_{014})]
\end{split}
\quad \raise1.4cm\hbox{\xymatrix{
& 0 \ar @{<-}[ddl] \ar @{<-}[ddrrr] \ar @{<-}[dr]^{g_1} \ar @{<-}[d] & &  \\
& 4\ar@{->}[dl]^{g_4} \ar@{->}[drrr]& 1 \ar @{<-}[dll] \ar @{<-}[drr]^{g_2} \ar @{<-}[l] & \\
3 \ar@{->}[rrrr]^{g_3}  & & & & 2 }}
\]
where $\eta_{234}=m_0 (\eta_{034}, \eta_{024}, \eta_{023})$ and $\eta_{124}=m_0(\eta_{024}, \eta_{014}, \eta_{012})$;
\[
\begin{split}
  &[((\eta_{012}, 1, 1), (\eta_{023}, 1), \eta_{034})]  \\
\overset{a\times id}{\mapsto}
&[((1, \eta_{123}, 1), (\eta_{013}, 1), \eta_{034}) ] \\
\overset{id \times a}{\mapsto}
&[((1, \eta_{123}, 1), (1, \eta_{134}), \eta_{014})] \\
\overset{a\times id} {\mapsto} &[((1,1, \eta_{234}), (1,
\eta_{124}), \eta_{014})]
\end{split}
\]
where $\eta_{123}= m_0(\eta_{023}, \eta_{013}, \eta_{012})$ and
$\eta_{134}= m_0(\eta_{034}, \eta_{014}, \eta_{013})$. Therefore,
the last map tells us that
\[ \eta_{123}= m_3 ( \eta_{234}, \eta_{134}, \eta_{124}). \]
Therefore associativity holds!

\subsubsection*{Comments on the \'etale condition} It is easy to see
that if $G_1\rra G_0$ is \'etale, by principality of the right $G$
action on $E_m$, the map $E_m\to G_0\times_M G_0$ is \'etale.
Moreover since $E_m \to \Lambda(X)_{2,j}=\Lambda[2,j](X)$ is a
surjective submersion by $Kan(2,j)$, by dimension counting, it is
furthermore an \'etale map.

\section{From Lie 2-groupoids to stacky Lie groupoids}\label{sec:2-slie}
If $X$ is a Lie 2-groupoid, then $G_1:= d_2^{-1}(s_0(X_0))\subset
X_2$, which is the set of bigons, is a Lie groupoid over
$G_0:=X_1$ (Lemma \ref{lemma:g1-g0}). Here we might notice that
there is another natural choice for the space of bigons, namely
$\tG_1:= d_0^{-1}(s_0(X_0))$. But $G_1 \cong \tG_1$ by the
following observation: given an element $\eta_3 \in G_1$, it fits
as the face opposite to 3 in diagram \eqref{pic:4gon} with $1\to
0$ and $2\to 3$ degenerate and $\eta_2$, $\eta_1$ degenerate; then
$m_0$ gives a morphism $\varphi: G_1 \to \tG_1$ and $m_3$ gives
the inverse. Therefore we might consider only $G_1$. Then $G_1
\rra G_0$ presents a stack which has an additional groupoid
structure. Hence we have the statement from 2-groupoids to stacky
Lie groupoids:

\begin{thm}\label{2-to-slie}
A Lie 2-groupoid (respectively 2-\'etale Lie 2-groupoid) $X$
corresponds to an SLie (respectively W-) groupoid $\cG\rra X_0$
where $\cG$ is presented by the Lie groupoid $G_1\rra G_0$.
\end{thm}

We prove this theorem by several lemmas.

\subsection{The stack $\cG$}
\begin{lemma}\label{lemma:g1-g0}
$G_1\rra G_0$ is a Lie  groupoid.
\end{lemma}
\begin{proof} The target and source maps are given by $d^1_0$ and
$d^1_1$. The identity $G_0\to G_1$ is given by $s^1_0: X_1 \to
X_2$. The image of $s^1_0$ is in $G_1 (\subset X_2)$. Their
compatibility conditions are implied by the compatibility
conditions of the structure maps of simplicial manifolds. This
will in particular imply that the identity is an embedding and
$\bs$ and $\bt$ are surjective submersions after we establish the
multiplication. The multiplication is given by the
3-multiplication of $X$ by the following picture:
\begin{equation}\label{pic:4gon}
\xymatrix{ &  & 0  & & \\
1 \ar[urr] & & & & 3 \ar[llll] \ar[llu] \ar[llld] \\
&  2 \ar[ruu] \ar[lu] & & & }
\end{equation}
where $\eta_3$ and its edges     $1\to 0$, $2\to 0$ and $2\to 1$
are degenerate (one should imagine them being very short).  Let
$\eta_i$ be the face facing $i$. Then any $(\eta_0, \eta_2)\in
G_1\times_{\bs, G_0, \bt} G_1$ fits in the above picture.  We
define $\eta_0\cdot \eta_2 = m_1(\eta_0, \eta_2, \eta_3)$ with
$\eta_3$ is the degenerate face in $s_0 \circ s_0(X_0)$
corresponding to the point 0(=1=2). Then the associativity of the
3-multiplications ensures the associativity of ``$\cdot$''. The
inverse map is also given by 3-multiplications: $\eta_2^{-1} =
m_0(\eta_1, \eta_2, \eta_3)$ with $\eta_1$ the degenerate face in
$s^1_0(X_1)$.
\end{proof}
\begin{remark}
Similar construction shows that $\tG_1\rra G_0$ with $\bt=d^1_2$,
$\bs=d^1_1$ is a Lie groupoid isomorphic to $G_1\rra G_0$ via the
map $\varphi^{-1}$.
\end{remark}

\subsection{Proof that $\cG \rra M$ is an SLie groupoid}
\subsubsection*{Source target maps and multiplication} There are
three maps $d^2_i: X_2\to X_1=G_0$ and they (as the moment maps of
the action) all correspond to a groupoid action respectively. The
actions are similarly given by the 3-multiplications as the
multiplication of $G_1$. The axioms of the actions are given by
the associativity. For example, for $d^2_1$, any $(\eta_0, \eta_2)
\in X_2 \times_{d^2_1, X_1, \bt_G}G_1$ fits inside picture
\eqref{pic:4gon}, but one imagines that only $1\to 0$ is a short
($\eta_3$ degenerate) edge.  Then
\begin{equation}\label{eq:g1-r-action}
\eta_0 \cdot \eta_2 := m_1(\eta_0, \eta_2, s_0d_2(\eta_0)).
\end{equation}
Moreover, notice that the four ways to compose source target and
face maps $G_1 \underset{\bt_G}{ \overset{\bs_G}{\rra}} G_0
\underset{d^1_1}{\overset{d^1_0}{\rra }}X_0$ only give two
different maps: $d^1_0\bs_G$ and $d^1_1 \bt_G$. They are
surjective submersions since $d_i$'s and $\bs_G$ $\bt_G$ are so
and they give the source and target maps $\bbs, \bbt: \cG\rra X_0$
where $\cG$ is the differentiable stack presented by $G_1\rra
G_0$. Therefore $\bbs$ and $\bbt$ are also surjective submersions
(Lemma 4.2 in \cite{tz}). We use these two maps to form the
product groupoid
\begin{equation}\label{eq:gg}
G_1\times_{d^1_0\bs_G, X_0, d^1_1 \bt_G} G_1
\rra  G_0\times_{d^1_0,X_0,d^1_1} G_0
\end{equation}
which presents the stack
$\cG\times_{\bbs, X_0, \bbt} \cG$. Then we have the following
lemma:

\begin{lemma}\label{lemma:mul}$(X_2, d^2_2\times d^2_0, d^2_1)$ is an
  H.S. bimodule from the groupoid in \eqref{eq:gg} to $G_1 \rra
G_0$.
\end{lemma}
\begin{proof}
By $Kan(2,1)$, $d^2_2\times d^2_1$ is a surjective submersion from
$X_2$ to $ G_0\times_{d^1_0,X_0,d^1_1} G_0$, so we only have to
show that the right action of $G_1\rra G_0$ on $X_2$ is free and
transitive. This is implied by $Kan(3,j)$ and $Kan(3,j)!$
respectively. \\
{\em transitivity}: any $(\eta_1, \eta_0)$ such that
$d^2_0(\eta_1)=d^2_0(\eta_0)$ and $d^2_2(\eta_0)=d^2_2(\eta_1)$
fits inside \eqref{pic:4gon} with $\eta_3$ a degenerate face (only
$1\to 0$ the short edge).
Then there exists $\eta_2:=m_2(\eta_0, \eta_1, \eta_3) \in G_1$, making $\eta_0\cdot \eta_2 =\eta_1$. \\
{\em freeness}: if $(\eta_0, \eta_2) \in X_2 \times_{d_1, X_1,
\bt}G_1$ satisfies $\eta_0\cdot \eta_2$ (=$m_1(\eta_0, \eta_2,
\eta_3))$ $ = \eta_0$. Then $\eta_2=m_2(\eta_0, \eta_0, \eta_3)$
and $\eta_3$ is degenerate. Thus $m_2 (\eta_0, \eta_0, \eta_3)$ =
$s^1_0(3\to 1)$ is a degenerate face. Therefore $\eta_2=1$.
\end{proof}

Therefore $X_2$ gives a morphism $m: \cG\times_{X_0}\cG \to \cG$.

\begin{lemma}With the source and target maps constructed above, $m$ is a multiplication of $\cG\rra X_0$.
\end{lemma}
\begin{proof} By construction, it is clear that $\bbt \circ m =
\bbt\circ pr_1$ and $\bbs \circ m = \bbs \circ pr_2$, where $pr_i:
\cG\times_{\bbs, X_0, \bbt} \cG \to \cG$ are the projections (see
the left picture below).
\begin{equation} \label{pic:mul}
\begin{xy}
  *\xybox{(0,0);<3mm,0mm>:<0mm,3mm>::
  ,0
  ,{\xylattice{-5}{0}{-4}{0}}}="S"
  ,{(-10,-10)*{\circ}}, {(-12, -12)*{_1}}
     ,{(0,0)*{\bullet}}, {(0, 2)*{^{0= \bbt m = \bbt pr_1}}}, {(10, -10)*{\bullet}},
     {(15, -8)*{^{2= \bbs m = \bbs pr_2}}},
     {(-10, -10) \ar@{->} (0,0)}, { (10, -10) \ar@{->} (-10,-10)} , { (10, -10) \ar@{->} (0,0)}
\end{xy}
\quad
\begin{xy}
*\xybox{(0,0);<3mm,0mm>:<0mm,3mm>::
  ,0
  ,{\xylattice{-5}{0}{-4}{0}}}="S"
  ,{(-10,-10)*{\bullet}}, {(-10, -12)*{_1}},{(-8,
  -14)*{\bullet}},  {(-8,
  -16)*{_{1'<1}}},
     ,{(0,0)*{\bullet}}, {(0, 2)*{^{0}}}, {(10, -10)*{\bullet}},
     {(10, -12)*{^{2}}}, {(15, -4)*{\bullet}}, {(17,-5)*{^3}},
     {(-10, -10) \ar@{->} (0,0)}, {(-10, -10) \ar@{.>} (-8,-14)},{ (10, -10) \ar@{->} (-10,-10)},
      { (10, -10) \ar@{.>} (-8,-14)} , { (10, -10) \ar@{->}
     (0,0)}, {(15, -4)\ar@{->}(10, -10)}, {(15, -4)
     \ar@{.>}(-8,-14)},{(15, -4)\ar@{->}(0, 0)},{(15, -4)\ar@{->}(-10,
     -10)}, {(-8,-14)\ar@{.>} (0,0)}
\end{xy}
\quad
\begin{xy}
*\xybox{(0,0);<3mm,0mm>:<0mm,3mm>::
  ,0
  ,{\xylattice{-5}{0}{-4}{0}}}="S"
  ,{(-10,-10)*{\bullet}}, {(-10, -12)*{_1}},
     ,{(0,0)*{\bullet}}, {(0, 2)*{^{0}}}
,{(1,-7)*{^{0'>0}}}
,{(0,-4)*{\bullet}}
, {(10, -10)*{\bullet}}
 , {(10, -12)*{^{2}}}
 ,    {(15, -4)*{\bullet}}
 ,{(17,-5)*{^3}}
 , {(-10, -10) \ar@{->} (0,0)},
     { (10, -10) \ar@{->}(-10,-10)},
      { (10, -10) \ar@{->} (0,0)},
     {(15, -4)\ar@{->} (10, -10)},
     {(15, -4)\ar@{->}(0, 0)},
     {(15, -4)\ar@{->}(-10,-10)},
{(0,-4)\ar@{.>}(0,0)},{(-10,-10)\ar@{.>}(0,-4)},{(10,-10)\ar@{.>}(0,-4)},{(15,-4)\ar@{.>} (0,-4)}
\end{xy}
\end{equation}
To show the associativity, we reverse the argument in Section
\ref{sec:statement}. There, we used the 2-morphism $a$ to
construct the 3-multiplications. Now we use the 3-multiplications
and their associativity to construct $a$. Given the two H.S.
bibundles presenting $m\circ (m\times id)$ and $m\circ (id\times
m)$ respectively, we want to construct a map $a$ as in
\eqref{eq:a}, where $E_m=X_2$ and $M=X_0$. Given any element in
$(X_2\times_{G_0} G_1) \times_{G_0\times_{X_0} G_0}X_2
/G_1\times_{X_0} G_1$, as in Section \ref{sec:statement}, we can
write it in the form of $[(\eta_3, 1,\eta_1 )]$, with $(\eta_1,
\eta_2, \eta_3)\in \hom(\Lambda[3, 0], X)$ for some $\eta_2$. Then
we define
\[a([(\eta_3, 1, \eta_1 )]):= [(1, m_0(\eta_1, \eta_2,
\eta_3), \eta_2 )].\] As before, we need to show that:1) the
R.H.S. is indeed in $ (G_1\times_{G_0} X_2)\times_{G_0\times_{X_0}
G_0}X_2 /\sim$; 2) the definition of $a$ does not depend on the
choice of $\eta_1$, $\eta_3$ and $\eta_2$; 3) $a$ is smooth. The
argument is the same. Here we only show 2) which is less obvious.
First of all, if we choose a different $\teta_2$, since $(\eta_1,
\eta_2, \eta_3) $ and $(\eta_1, \teta_2, \eta_3) $ are both in
$\hom(\Lambda[3,0], X)$, we have $d^2_2(\eta_2)=d^2_2(\teta_2)$
and $d^2_1(\eta_2)=d^2_1(\teta_2)$. So $\eta_2=\eta_{013}$ and
$\teta_2=\eta_{0 1' 3}$ form a degenerate horn (see the second
picture of \eqref{pic:mul}).  By $Kan(3,0)$ there exists
$\gamma=\gamma_{1'13}\in G_1$ such that $(1, \gamma)\cdot
\teta_2(=\gamma \cdot \teta_2)=\eta_2$, that is
$\eta_{013}=m_1(\gamma, \eta_{01'3}, s^1_1(0\to 1))$. Then by
associativity and the definition of the $G_1$ action
\eqref{eq:g1-r-action}, we have $m_0(\eta_{023}, \eta_{01'3},
\eta_{01'2}) = \eta_{1'23}=\eta_{123} \cdot \gamma$. Therefore we
have $[(1, \eta_{1'23}, \teta_2)]=[(1, \eta_{123}, \eta_2)]$. So
the choice of $\eta_2$ will not affect the definition of $a$. If
we choose a different $(\teta_3=\eta_{0'12},1,
\teta_1=\eta_{0'23})$, then we can suppose
$\eta_3=\teta_3\gamma_{00'2}$ and
$\teta_1=(\gamma_{00'2},1)\cdot\eta_1=\gamma\cdot \eta_1$ (see the
third picture of \eqref{pic:mul}). Then $(\teta_1, \eta_2,
\teta_3) \in \hom(\Lambda[3,0], X)$ and
\[m_0(\teta_1, \eta_2, \teta_3) =\eta_{123}=m_0(\eta_1,\eta_2,\eta_3).\]
So this choice will not affect $a$ neither.

Now the higher coherence of $a$ follows from the associativity by
the same argument as in Section \ref{sec:3-asso}.
\end{proof}

\subsubsection*{Identity} Now we notice that $s_0: X_0
\hookrightarrow G_0$ and $e_G
\circ s_0: X_0 \hookrightarrow G_1$ with $e_G$  the identity of $G$ form a groupoid morphism from
$X_0\rra X_0$ to $G_1 \rra G_0$. It gives a morphism $\bar{e}:
X_0\to \cG$ on the level of stacks.
\begin{lemma} $\bar{e}$ is the identity of $\cG$.
\end{lemma}
\begin{proof} Recall from Definition \ref{def:sliegpd} that we need to show that there is a 2-morphism $b_l$
between the two maps $m\circ (\bar{e}\times id)$ and $pr_2$: $X_0
\times_{X_0, \bt} \cG \to \cG$, and similarly a 2-morphism $b_r$.
In our case, the H.S. bibundles presenting these two maps are
$X_2|_{X_0\times_{X_0} G_0}$ and $G_1$ respectively and they are
the same by construction, hence $b_l=id$. For $b_r$, we notice
that $X_2|_{G_0\times_{X_0}X_0} = \tG_1$ and the isomorphism $
\varphi^{-1}: \tG_1 \to G_1$ is $b_r$.

Now we need to show the higher coherences that $b$'s satisfy. The
proofs are similar and here we prove \ref{itm:br} in Definition
\ref{def:sliegpd} in detail. Translating the diagram
\eqref{diag:br} in the language of groupoids and bibundles, we
obtain
\[\small{\xymatrix{
G_1\times_{X_0}G_1\times_{X_0}X_0 \ar[dd] \ar@<-1ex>[dd] & &
G_1\times_{X_0}X_0 \ar[dd] \ar@<-1ex>[dd] & & G_1 \ar[dd]
\ar@<-1ex>[dd]
\\
& X_2\times_{X_0}X_0\ar[dl]\ar[dr] & & \tG_1 \ar[dl]\ar[dr] & \\
G_0\times_{X_0}G_0\times_{X_0}X_0 & & G_0\times_{X_0} X_0 & & G_0
\\
&  & G_0\times_{X_0}G_0 & & \\& G_1 \times_{X_0} \tG_1 \ar[uul]
\ar[ur] & & X_2 \ar[ul]
\ar[uur]\\
& & G_1\times_{X_0}G_1\ar[uu] \ar@<-1ex>[uu] &&}.}
\]
Corresponding to \eqref{eq:br}, we need to show that the following
diagram commute:
$$ \small{\xymatrix{ (X_2\mathop\times\limits_{X_0} X_0)
\mathop\times\limits_{G_0\mathop\times\limits_{X_0} X_0}
\tG_1/\sim \ar[r]^a \ar[d]_{b_r} & (G_1\mathop\times\limits_{X_0}
\tG_1 )\mathop\times\limits_{G_0\mathop\times\limits_{X_0} G_0}
X_2 /\sim \ar[r]^{b_r}
 & (G_1\mathop\times\limits_{X_0} G_1) \mathop\times\limits_{G_0 \mathop\times\limits_{X_0} G_0} X_2 /
\sim \ar[ld]^{\cong} \\
(X_2\mathop\times\limits_{X_0} X_0
)\mathop\times\limits_{G_0\mathop\times\limits_{X_0}X_0}G_1/\sim
\ar[r]^{\cong} & X_2 & }} $$
\[
 \xymatrix{ [(\eta_3, 1, \eta_1)]
\ar@{|->}[r] \ar@{|->}[d]& [(1,  \eta_0=1, \eta_2)]
\ar@{|->}[r] & [(1, \eta'_0=1, \eta_2 )] \ar[ld]^? \\
[(\eta_{012},1,\eta_{0'02})] \ar@{|->}[r] &
\eta_{012}\cdot\eta^{-1}_{0'02} }.\]
\begin{equation}\label{pic:br}
\begin{xy}
*\xybox{(0,0);<3mm,0mm>:<0mm,3mm>::
  ,0
  ,{\xylattice{-5}{0}{-4}{0}}}="S"
  ,{(-10,-10)*{\bullet}}, {(-10, -12)*{_1}}
     ,{(0,0)*{\bullet}}, {(0, 2)*{^{0}}}, {(10, -10)*{\bullet}},
     {(10, -12)*{^{2}}}, {(15, -8)*{\bullet}}, {(5, 0)*{\bullet}},
     {(7,2)*{^{0'<0}}}, {(17,-8)*{^3}},
     {(-10, -10) \ar@{->} (0,0)}, {(-10, -10) \ar@{->} (5,0)},{ (10, -10) \ar@{->} (-10,-10)} , { (10, -10) \ar@{->} (5,0)} , { (10, -10) \ar@{->}
     (0,0)}, {(15, -8)\ar@{.>}(10, -10)}, {(15, -8)
     \ar@{.>}(5,0)},{(15, -8)\ar@{.>}(0, 0)},{(15, -8)\ar@{.>}(-10,
     -10)}, {(5,0)\ar@{<-} (0,0)}
\end{xy}
\end{equation}
Let us explain the diagram: A $[(\eta_3, 1, \eta_1)] \in
(X_2\times_{X_0} X_0) \times_{G_0\times_{X_0} X_0} \tG_1/\sim$
fits inside \eqref{pic:br} with $3\to 2$ a degenerate edge and
$\eta_0=\eta_{123}$ a degenerate face\footnote{Now when it is
confusing, we call a face by its three vertices, for example now
$\eta_1=\eta_{023}$.}, therefore its image under $a$ is $[(1,
\eta_0=1, \eta_2)]$. By the construction of $b_r$,
$\eta'_0=b_r(\eta_0)=1$ is also degenerate. Then $[(1, \eta'_0=1,
\eta_2)]$ maps to $\eta_2$ under ``$\cong$''. Now we follow the
other direction of the maps. Then $b_r(\eta_{023})=\eta_{0'02}$.
The image of $[(\eta_{012},1,\eta_{0'02})]$ under ``$\cong$'' is
$\eta_{012}\cdot\eta_{0'02}^{-1} = m_1(\eta_{0'02}^{-1},
\eta_{0'01}=1, \eta_{012}) = \eta_{0'12}$. Notice that most of the
``flat'' faces are degenerate except for $\eta_{0'02}$ and
$\eta_{023}$. Therefore consider the 3-simplices $(0',0,1,3)$ and
$(0',1,2,3)$, then we have $\eta_{013}=\eta_{0'13}$ and
$\eta_{0'13}=\eta_{0'12}$ by \eqref{coco}. Since $\eta_{013}$ is $
\eta_2$, the desired diagram commutes.
\end{proof}

\subsubsection*{Inverse} By Prop. \ref{prop:inverse}, we only have
to show that the actions of $G$ and $G^{op}$ on $X_2
\times_{d^2_1, G_0} M$, induced respectively by the first and
second components of the left action of $G_1 \times_M G_1 \rra
G_0\times_M G_0$, are principal (see \eqref{eq:inv-construct}). We
prove this for the first copy of $G_1$ and the proof for the
second is similar. A $(\eta_3, \eta_{1}) \in G_1 \times_{\bs_G,
G_0,d^2_0} X_2 \times_{d^2_1, G_0} M $ fits inside
\eqref{pic:4gon} with $\eta_2$ a degenerate face corresponding to
the point $0=1=3$. Then the freeness of the action is implied by
$Kan!(3, 0)$ and the transitivity of the action is implied by
$Kan(3, 0)$.

\subsubsection*{Comments on the \'etale condition} If $X_2 \to
\Lambda[2,j](X)$ are \'etale maps, then the moment map $J_l:
E_m\cong X_2 \to G_0\times_{\bs,M, \bt} G_0$ is \'etale. By
principality of the right action of $G_1 \rra G_0$, it is an
\'etale groupoid. This concludes the proof of the Theorem
\ref{2-to-slie}.

\section{Equivalences of Lie 2-groupoids}

As in Theorem \ref{thm:slie-2}, from an SLie groupoid to a Lie
2-groupoid, one has to choose a chart of the stack $\cG$.
Therefore to construct a 1-1 correspondence between SLie groupoids
and Lie 2-groupoids, we have to quotient by some sort of
equivalence of 2-groupoids since one could have different
2-groupoids $X$ and $Y$ corresponding to different charts $ G_0$
and $H_0$ of the same stack $\cG$. In this section, we establish
various equivalences of Lie 2-groupoids.

We first look at the more general case of \lht s and motivate our
later definitions. The proof of this direction itself is carried
out in a later work \cite{hz2}.

\subsection{Strict maps of \lht s}\label{sec:equi-lht}

The reader's first guess is probably that a morphism $f:X\to Y$ of
\lht s ought to be a simplicial smooth map i.e. a collection of
smooth maps $f_n:X_n\to Y_n$ that commute with faces and
degeneracies. In the language of categories, this is just a
natural transformation from the functor $X$ to the functor $Y$. We
shall call such a natural transformation a \em strict map \rm from
$X$ to $Y$. Unfortunately, it is known that, already in the case
of usual Lie groupoids, such strict notions are not good enough.
Indeed there are strict maps that are not invertible even though
they ought to be isomorphisms. That's why people introduced the
notion of {\em Hilsum-Skandalis maps} \cite{hs}. Here is an
example of such a map: consider a manifold $M$ with an open cover
$\{\U_\alpha\}$. The simplicial manifold $X$ with
$X_n=\bigsqcup_{\alpha_1,\ldots,\alpha_n}\U_{\alpha_1}\cap\cdots\cap\U_{\alpha_n}$
maps naturally to the constant simplicial manifold $M$. All the
fibers of that map are simplices, in particular they are
contractible simplicial sets. Nevertheless, that map has no
inverse.

The second guess is then to define a special class of strict maps
which we shall call {\em equivalences}. A map from $X$ to $Y$
would then be a {\em zig-zag}  of strict maps
$X\stackrel{\sim}{\leftarrow}Z\to Y$, where the map $Z\to X$ is
one of these equivalences.

This will not be our final choice for what a morphism of \lht s
should be. The notion of equivalence is nevertheless very useful
(e.g. to define sheaf cohomology of \lht s) and we will study it
in this section. It is very much inspired from the notion of
equivalence of simplicial sets (sometimes called weak
equivalence). Recall (\cite{goerss-jardine}) that a map $S\to T$
of simplicial sets is an equivalence if it induces isomorphisms of
all the topological homotopy groups. Here is an equivalent
condition which we will generalize to Lie homotopy types:

\begin{lemma}\label{eqss}
A map $S\to T$ of Kan simplicial sets is an equivalence if, for
any $m\ge 0$ and any commutative solid arrow diagram
\begin{equation}\label{eqeq}
\xymatrix{
\partial\D^n \ar[r] \ar@{^{(}->}[d] &S \ar[d]\\
\D^n \ar[r] \ar@{.>}[ru] & T}
\end{equation}
there exists a dotted arrow that makes both triangles commute.
Here $\partial \D^n$ stands for the boundary of the $n$-simplex.
\end{lemma}

\begin{proof}
Let $\alpha\in\pi_{n-1}(S)$ be represented by some map
$\partial\D^n \to S$ and assume $\alpha\mapsto 0\in\pi_{n-1}(T)$.
This means that we have a map $\D^n\to T$ that makes diagram
(\ref{eqeq}) commute. By hypothesis,  we get a map $\D^n\to S$
therefore $\alpha=0\in\pi_{n-1}(S)$. This proves the injectivity
of $\pi_{n-1}(S)\to\pi_{n-1}(T)$.

Now let us consider an element $\beta\in\pi_n(T)$ represented by a
map $(\D^n,\pD^n)\to(T,*)$. That map fits into a diagram
(\ref{eqeq}) where the top arrow sends everything to the base
point. By hypothesis we get a map $(\D^n,\pD^n)\to(S,*)$. It
represents an element of $\pi_n(S)$ that is mapped to $\beta$.
This proves the surjectivity of $\pi_n(S)\to\pi_n(T)$.

The Kan condition has been used implicitly to find representatives
of the various homotopy classes (otherwise we would have needed to
subdivide the simplices).
\end{proof}

Translating the condition of lemma \ref{eqss} into hom spaces
gives:

\begin{df}\label{defequivalence}
A strict map $f:Z\to X$ of \lht s is an equivalence if the natural
maps from $Z_n=\hom(\D^n,Z)$ to the pull-back
$PB\big(\hom(\pD^n,Z)\to\hom(\pD^n,X)\leftarrow\hom(\D[n],X)\big)$
are surjective submersions for all $n\ge 0$.
\end{df}

The proof of this being an equivalence is in \cite{hz2}. But even
before this, we will need to talk a lot about spaces of the form
$PB\big(\hom(A,Z)\to\hom(A,X)\leftarrow\hom(B,X)\big)$, where the
maps are induced by some fixed maps $A\to B$ and $Z\to X$. To
avoid the cumbersome pull-back notation, we shall denote these
spaces by
$$\PB(A,Z,B,X).$$
This notation indicates that the space parameterizes all commuting
diagrams of the form
$$\begin{matrix}A&\!\!\!\longrightarrow&\!\!\!Z\\
\downarrow&&\!\!\!\downarrow\\
B&\!\!\!\longrightarrow&\!\!\!X,\end{matrix}$$ where we allow the
horizontal arrows to vary but we fix the vertical ones. 

As in the case of definition \ref{defngroupoids}, we need to
justify why the pull-backs used in definition \ref{defequivalence}
are manifolds. This is specially surprising since the spaces
$\hom(\pD^n,Z)$ need not be manifolds (for example take $n=2$ and
$Z$ the cross product groupoid associated to the action of $S^1$
on $\R^2$ by rotation around the origin).

\begin{lemma}\label{lemma:indct-mfd-equi}
Let $S$ be a finite collapsible simplicial set and $T\to S$ a
sub-simplicial set of dimension $\le k$. Let $f:Z\to X$ be a
strict map of \lht s such that the natural map
$$
Z_n\rightarrow\PB(\pD^n,Z,\D^n,X)
$$
is a surjective submersion for all $n\le k$ (assuming the right
hand sides are known to be manifolds). Then the space
$$\PB(T,Z,S,X)$$is naturally a manifold.
\end{lemma}

\begin{proof}
Let $T'$ be a sub-simplicial set obtained by deleting one simplex
from $T$ (without its boundary). We have a push-out diagram

\begin{diagram}
T'      &    \rTo    & \SWpbk T\\
 \uTo       &   & \uTo \\
\pD^n&    \rinto    &\D^n&.
\end{diagram}

Similarly to the proof of lemma \ref{collapsable1} this gives a
pull-back diagram
\begin{diagram}
\PB(T',Z,S,X)      &&   \lTo    & \raisebox{-.7cm}{\SWpbk} \PB(T,Z,S,X)&\phantom{=Z_n}\\ \\
 \dTo       &   && \dTo \\
\PB(\pD^n,Z,S,X)&&    \lTo    &\PB(\D^n,Z,S,X)&,
\end{diagram}

which may be combined with the pull-back diagram

\begin{diagram}
\PB(\pD^n,Z,S,X)      &&   \lTo    & \raisebox{-.7cm}{\SWpbk} \PB(\D^n,Z,S,X)\\ \\
 \dTo       &   && \dTo \\
\PB(\pD^n,Z,\D^n,X)&&    \lTo    &\PB(\D^n,Z,\D^n,X)&=Z_n
\end{diagram}

to give yet another pull-back diagram

\begin{equation}\label{bigpb}
\begin{diagram}
\PB(T',Z,S,X)      &&   \lTo    & \raisebox{-.7cm}{\SWpbk} \PB(T,Z,S,X)\\ \\
 \dTo       &   && \dTo \\
\PB(\pD^n,Z,\D^n,X)&&    \lTo    &\PB(\D^n,Z,\D^n,X)&=Z_n.
\end{diagram}
\end{equation}

By induction on the size of $T$, and by lemma \ref{collapsable1}
(case $T=\emptyset$), we may assume that the upper left and lower
left spaces in (\ref{bigpb}) are known to be manifolds. The bottom
arrow is a submersion by hypothesis. Therefore by transversality,
 the upper right space is also a manifold, which is what we wanted to prove.
\end{proof}

\subsection{Morita Equivalence of Lie 2-groupoids}
Now we define equivalences of Lie 2-groupoids based on Section
\ref{sec:equi-lht}. To simplify our notation, $\thra$ and $\thla$
always denote surjective submersions.

\begin{defi}\label{defi:equi-2gpd}
Given two Lie 2-groupoids $Z_2 \Rrightarrow Z_1 \rra Z_0$ and $X_2
\Rrightarrow X_1 \rra X_0$, a strict morphism $f$ between them is
an equivalence if the natural map from $Z_n$ to the pull-back
$PB(\hom(
\partial \Delta^n , Z)\ra \hom(\partial \Delta^n, X) \la X_n)$
are surjective submersions for all $n= 0,1$ and an isomorphism for
$n=2$.
\end{defi}
As in the case of Lie homotopy types, we find that $PB(\hom(
\partial \Delta^n , Z)\ra \hom(\partial \Delta^n, X) \la X_n)$ is
a manifold.

Similarly to the pull-back Lie groupoid, we make the following
observation: Let $X$ be a 2-groupoid and $Z_1\rra Z_0$ be two
manifolds with structure maps as in \eqref{eq:face-degen} up to
the level $n\leq 1$, and $f_n: Z_n \to X_n$ preserving the
structure maps $d^{n}_k$'s and $s^{n-1}_k$'s for $n\leq 1$. Then
$\hom(
\partial \Delta^n , Z)$ still makes sense for $n\leq 1$. We further suppose
$f_0: Z_0 \twoheadrightarrow X_0$ (hence $Z_0\times Z_0
\times_{X_0\times X_0} X_1$ is a manifold) and $Z_1\thra Z_0\times
Z_0 \times_{X_0\times X_0} X_1$. That is to say  that the induced
map from $Z_n$ to the pull-back $PB(\hom(
\partial \Delta^n , Z)\ra \hom(\partial \Delta^n, X) \la X_n)$
are surjective submersions for $n= 0, 1$. Then we form
\[Z_2= PB(\hom(
\partial \Delta^2 , Z)\ra \hom(\partial \Delta^2, X) \la X_2).\]
It is easy to see that the proof of Lemma
\ref{lemma:indct-mfd-equi} still guarantees that $Z_2$ a manifold.
Moreover there are $d^2_i: Z_2\to Z_1$ induced by the natural
projections $\hom(\partial \Delta^2 , Z)\to Z_1$; $s^1_i: Z_1 \to
Z_2$ by
\[ s^1_0(h)=(h,h,s^0_0(d^1_0(h)),s^1_0(f_1(h))), \quad
s^1_1(h)=(s^0_0(d^1_1(h)),h,h,s^1_1(f_1(h)));
\]
$m_i: \Lambda(Z)_{3, i}\to Z_2$ by for example
\[ m_0(( h_2, h_5, h_3, \bareta_1), (h_4, h_5, h_0, \bareta_2), (h_1, h_3, h_0, \bareta_3))=
( h_2, h_4, h_1, m_0(\bareta_1, \bareta_2, \bareta_3)), \] and
similarly for other $m$'s.
\[ \xymatrix{ &  & 0  & & \\
1 \ar[urr]^{h_0} & & & & 3 \ar[llll]^{h_4} \ar[llu]^{h_5} \ar[llld]^{h_2} \\
&  2 \ar[ruu]^{h_3} \ar[lu]^{h_1} & & & } \] Then $Z_2
\Rrightarrow Z_1 \rra Z_0$ is a Lie 2-groupoid and we call it the
{\em pull-back 2-groupoid} by $f$. Moreover $f: Z\to X$ is an
equivalence with the natural projection $f_2: Z_2 \to X_2$.

Definition \ref{defi:equi-2gpd} can be easily generalized to
higher dimensions.
\begin{defi}\label{defi:equi-ngpd}
A strict map $f: Z \to X$ of Lie $m$-groupoids is an equivalence
if the natural maps from $Z_n$ to the pull-back $PB(\hom(
\partial \Delta^n , Z)\ra \hom(\partial \Delta^n, X) \la X_n)$
are surjective submersions for all $ 0\leq n \leq m-1 $ and
isomorphism for $n= m$.
\end{defi}

\begin{remark}
First of all, it is easy to see that when $m=1$, we get the
equivalence (or pull-back) of Lie groupoids. Moreover, since for a
Lie $m$-groupoid layers higher or equal to $m+1$ are determined by
layers lower or equal to $m$, the maps from $Z_n$ to the pull-back
$PB(\hom(
\partial \Delta^n , Z)\ra \hom(\partial \Delta^n, X) \la X_n)$ are
isomorphisms for all $n\geq m$.
\end{remark}

Similar to the proof that $PB(\hom(
\partial \Delta^n , Z)\ra \hom(\partial \Delta^n, X) \la X_n)$ is
a manifold, we have the following lemma:

\begin{lemma} \label{lemma:proj-sub}The projection $PB(\hom(
\partial \Delta^n , Z)\ra \hom(\partial \Delta^n, X) \la X_n) \to
X_n$ is a surjective submersion.
\end{lemma}
\begin{proof} We can use the same induction as in Lemma \ref{lemma:indct-mfd-equi}  and only have to
notice that the upper lever map in \eqref{bigpb} is a surjective
submersion (then let $T'=\emptyset$).
\end{proof}

\begin{lemma} \label{lemma:comp-equi}
 The composition of equivalences is still an
equivalence.
\end{lemma}
\begin{proof}
Suppose that $Z \overset{f}{\to} Y \overset{g}{\to} X$ are two
equivalences. Then we have
\[
\xymatrix@C=0.3cm{ Z_n \ar@{->>}[r]  & \hom(\partial \Delta^n,
Z)\times_{\hom(\partial \Delta^n, Y)} Y_n
\ar[r]\ar@{->>}[d]^{\Leftarrow}
& Y_n \ar@{->>}[d] & \\
& {\hom(\partial \Delta^n,Z)\times_{\hom(\partial \Delta^n, X)}
X_n} \ar[d]\ar[r] & {\hom(\partial
\Delta^n,Y)\times_{\hom(\partial \Delta^n, X)} X_n}
\ar[d]\ar[r] & {X_n} \ar[d] \\
& {\hom(\partial \Delta^n,Z)} \ar[r] & {\hom(\partial \Delta^n,Y)}
\ar[r] & {\hom(\partial \Delta^n, X)} }
\]
Therefore we have $Z_n \thra \hom(\partial
\Delta^n,Z)\times_{\hom(\partial \Delta^n, X)} X_n$ when $n\leq
1$. When $n=2$ one just has to replace $\thra$ by $\cong$ in the
above diagram and the result. 
\end{proof}
\begin{remark}
It is not hard to see this result generalizes to $n$-groupoids.
\end{remark}

\begin{lemma}\label{lemma:inverse-comp-equi}
If $f: Z \to X$ and $f': Z \to X$ are equivalences, then there
exist a Lie 2-groupoid $Z''=Z \times_{X} Z'$ and
two equivalences $Z \lla Z'' \lra Z'$.
\end{lemma}
\begin{proof}
First of all we prove that $Z''_n = Z_n \times_{X_n} Z'_n$ (a
manifold by Lemma \ref{lemma:proj-sub}) with product face and
degeneracy maps is a Lie 2-groupoid. We use Prop-Def.
\ref{def:finite-2gpd}. Notice that the maps on both sides of  $
Z_0 \lla Z_0\times_{X_0} X_1 =PB(Z_0 \thra X_0 \thla X_1) \lra
X_1$ are surjective submersions by the property of pull-backs.
Similarly we have a surjective submersion,
\[Z_0 \times Z_0 \times_{X_0 \times X_0} X_1 = PB(
(Z_0\times_{X_0}X_1) \thra X_1 \thla (Z_0\times_{X_0}X_1))
\twoheadrightarrow (Z_0\times_{X_0} X_1),
\] hence $Z_1 \thra Z_0\times_{X_0} X_1$ by composing with $Z_1 \thra Z_0 \times Z_0 \times_{X_0 \times X_0} X_1 $.
Similarly we have this for $Z'$. For the 1-Kan conditions of
$Z''$, we should show that the degeneracy maps $Z_1\times_{X_1}
Z'_1\to Z''_0=Z_0\times_{X_0} Z'_0 $ are surjective submersions.
The degeneracy map $d^1_0 \times d^1_0$ is a surjective submersion
by chasing the arrows denoted with ``!'' in the following diagram,
\[
\xymatrix@C=0.3cm@R=0.3cm{{Z_1\times_{X_1} Z'_1}\ar[rr]
\ar@{->>}[rd]^{\Leftarrow}_{!} & & {*}\ar[rr]
\ar@{->>}[rd]^{\Leftarrow}  & & {Z'_1} \ar@{->>}[rd] & &
\\
& {*} \ar@{->>}[rr]_{\Uparrow}^{!} \ar[rd] & & {*} \ar@{->>}[rr]
\ar[rd] \ar@{->>}'[d]_{!}^{\Leftarrow}[dd]
& & {Z'_0 \times_{X_0} X_1}\ar@{->>}[rd] \ar@{->>}'[d][dd] &  \\
& & {Z_1} \ar@{->>}[rr] & & {Z_0 \times_{X_0}X_1} \ar@{->>}[rr]
\ar@{->>}^{\Leftarrow}[dd]
& & {X_1} \ar@{->>}[dd]^{d_0} \\
& & & {Z_0\times_{X_0} Z'_0} \ar'[r][rr] \ar[rd] & & {Z'_0} \ar[rd] & \\
& & & & {Z_0} \ar[rr]  & & {X_0} }
\] where $*$ denotes non-important pull-backs.
The same diagram proves the result for $d^1_1\times d^1_1$ also.
This argument in fact shows the following general result about
surjective submersions of manifolds,
\begin{claim} \label{claim:ss-tool} If $U\thra \bar{U}$, and the map $V\to \bar{V}$ factors through a surjective
submersion $V\thra \bar{V} \times_{\bar{U}} U$, and $W\to \bar{W}$
factors through $W \thra U\times_{\bar{U}}\bar{W}$, then the
natural map $V\times_{U} W \thra \bar{V}\times_{\bar{U}}\bar{W}$
is a surjective submersion. Here we suppose all the fibre products
to be manifolds.
\end{claim}

For the 2-Kan condition of $Z''$, we need to show that the maps
$Z''_2 \to \hom(\Lambda[2,j], Z'')$ are surjective submersions.
Since $\hom(\Lambda[2,j], Z'')=\hom(\Lambda[2,j],
Z)\times_{\hom(\Lambda[2,j], X)} \hom(\Lambda[2,j], Z')$, by Claim
\ref{claim:ss-tool} and $X_2 \thra \hom(\Lambda[2,j], X)$, we only
have to show $Z_2 \thra \hom(\Lambda[2,j], Z)
\times_{\hom(\Lambda[2,j], X)} X_2$ and the similar statement for
$Z'$. Notice that
\[ Z_0\times Z_0
\thra X_0\times X_0 ,\quad Z_1 \thra Z_0\times
Z_0\times_{X_0\times X_0} X_1, \] by assumption. Therefore
\[
\begin{split}
 \hom(\partial \Delta^2, Z) & =\hom(\Lambda[2,j],
Z)\times_{Z_0\times Z_0} Z_1 \\
&\thra \hom(\Lambda[2,j], Z)\times_{Z_0\times Z_0}
Z_0\times Z_0\times_{X_0\times X_0} X_1 \\
&=\hom(\Lambda[2,j], Z)\times_{X_0\times X_0} X_1 \\
&=\hom(\Lambda[2,j], Z) \times_{\hom(\Lambda[2,j], X)}
\hom(\partial \Delta^2, X),
\end{split}
\]
where in the last step we use again $\hom(\partial \Delta^2,
X)=\hom(\Lambda[2,j], X)\times_{X_0\times X_0} X_1$. Hence
\[
\begin{split}
Z_2 & = \hom(\partial \Delta^2, Z)
\times_{\hom(\partial \Delta^2, X)} X_2 \\
&\thra \hom(\Lambda[2,j], Z) \times_{\hom(\Lambda[2,j], X)}
\hom(\partial \Delta^2, X)
\times_{\hom(\partial \Delta^2, X)} X_2 \\
&=\hom(\Lambda[2,j], Z) \times_{\hom(\Lambda[2,j], X)} X_2.
\end{split}
\]
The 3-multiplications on $Z''_2$ are induced by those of $X_2$, by
the second row of the equation,
\begin{equation}\label{eq:z''-2}
\begin{split}
Z''_2 & = \hom(\partial \Delta^2, Z) \times_{\hom(\partial
\Delta^2, X)} X_2 \times_{X_2} \hom(\partial\Delta^2, Z')
\times_{\hom(\partial \Delta^2, X)} X_2 \\
 & = \hom(\partial \Delta^2, Z) \times_{\hom(\partial
\Delta^2, X)} X_2\times_{\hom(\partial \Delta^2, X)}
\hom(\partial\Delta^2, Z') \\
& = Z_2 \times_{\hom(\partial \Delta^2, X)} \hom(\partial\Delta^2,
Z') \\
& = \hom(\partial \Delta^2, Z'') \times_{\hom(\partial \Delta^2,
Z)} Z_2.\\
\end{split}
\end{equation}
In the last row, we use $\hom(\partial \Delta^2, Z'')
=\hom(\partial \Delta^2, Z) \times_{\hom(\partial \Delta^2, X)}
\hom(\partial\Delta^2, Z')$. Therefore, the 3-multiplications of
$Z''_2$ satisfy all the axioms as they are essentially the
3-multiplications of $X_2$. Now we have proven that $Z''$ is a Lie
2-groupoid. Moreover, from the last row, we can also see that the
projection $Z''_2 \to Z_2$ preserves the 3-multiplications since
both of theirs are induced by those of $X_2$. An analogous result
holds for $Z'$. Therefore to complete the proof, we only have to
show \[
\begin{split}
Z''_1 = Z_1\times_{X_1} Z'_1 & \thra \hom(\partial \Delta^1,
Z'')\times_{\hom(\partial \Delta^1, Z)} Z_1  \\
& = \hom(\partial \Delta^1, Z')\times_{\hom(\partial
\Delta^1, X)} Z_1,
\end{split}
\]since $\hom(\partial \Delta^1, Z'')
=\hom(\partial \Delta^1, Z) \times_{\hom(\partial \Delta^1, X)}
\hom(\partial\Delta^1, Z')$. But this can be easily deduced from
the following pull-back diagram,
\[
\xymatrix{Z_1\times_{X_1} Z'_1 \ar[r] \ar@{->>}[d]^{\Leftarrow}&
Z'_1\ar@{->>}[d] \\
Z_1\times_{ \hom(\partial \Delta^1, X)}\hom(\partial \Delta^1,
Z') \ar[r] \ar[d]& \hom(\partial \Delta^1,
Z')\times_{\hom(\partial \Delta^1, X)} X_1 \ar[d] \\
Z_1 \ar[r]& X_1,}
\] and the fact that $Z'_1\thra
\hom(\partial \Delta^1, Z')\times_{\hom(\partial \Delta^1, X)}
X_1$.
\end{proof}

\begin{defi}\label{defi:m-equi-2gpd} Two Lie 2-groupoids $X$ and $Y$ are
{\em Morita equivalent} if there is a Lie 2-groupoid $Z$ such that
both of the maps $X \lla Z \lra Y$ are equivalences. With the
above two lemmas, this definition does give an equivalence
relation. We call it {\em Morita equivalence} of Lie 2-groupoids.
\end{defi}

\subsection{1-Morita Equivalence of Lie 2-groupoids}
However Morita equivalent Lie 2-groupoids correspond to Morita
equivalent SLie groupoids \cite{bz}. Hence to obtain isomorphic
SLie groupoids, we need a stricter notion of equivalence of Lie
2-groupoids.

\begin{pdef}\label{defi:1t-m-equi-2gpd}
A strict map of 2-groupoids $f: Z \to X$ is a {\em 1-equivalence}
if it is an equivalence with $f_0$ an
isomorphism. 
Two 2-groupoids $X$ and $Y$ are 1-Morita equivalent if there is a
2-groupoid $Z$ such that both of the maps $X \lla Z \lra Y$ are
1-equivalences. It gives an equivalence relation between Lie
2-groupoids and we call it {\bf 1-Morita equivalence}.
\end{pdef}
\begin{proof}
It is easy to see that the composition of 1-equivalences is still
a 1-equivalence.  We just have to notice that in Lemma
\ref{lemma:inverse-comp-equi} if both equivalences are
1-equivalences, then the projections $Z_0\la Z_0\times_{X_0} Z'_0
\ra Z'_0$ are isomorphisms since $Z_0\times_{X_0} Z'_0\cong
Z_0\cong X_0 \cong Z'_0$.
\end{proof}
\begin{remark} \label{rk:1t-m-equi}
For a 1-equivalence $Z\to X$, since $f_0: Z_0 \cong X_0$, we have
$\hom(\partial \Delta^1, Z) = \hom(\partial \Delta^1, X)$. So the
condition on $f_1$ in Definition \ref{defi:equi-2gpd} becomes
$f_1: Z_1 \thra X_1$.
\end{remark}

\section{One-to-one correspondence between Lie 2-groupoids and SLie groupoids}

In this section, we use two lemmas to prove the following theorem:
\begin{thm}\label{thm:1-1} There is a 1-1 correspondence of Lie 2-groupoids
(respectively 2-\'etale Lie 2-groupoids) modulo 1-Morita
equivalence and SLie groupoids (respectively W-groupoids).
\end{thm}

W-groupoids are isomorphic if and only if they are isomorphic as
SLie groupoids, and 1-Morita equivalent 2-\'etale Lie 2-groupoids
are 1-Morita equivalent Lie 2-groupoids. Therefore the \'etale
version of the theorem is implied by the general case and we only
have to prove the general case.

For the lemma below, we fix our notation: $X$ and $Y$ are Lie
2-groupoids in the sense of Prop-Def. \ref{def:finite-2gpd};
$G_0=X_1$ and $H_0=Y_1$; $X_0=Y_0=M$. $G_1$ and $H_1$ are the
spaces of bigons in $X$ and $Y$ respectively, that is
$d_2^{-1}(s_0(Y_0))$; Both $G_1 \rra G_0$ and $H_1 \rra H_0$ are
Lie groupoids and they present differentiable stacks $\cG$ and
$\cH$ respectively. Moreover $\cG \rra M$ and $\cH \rra M$ are
stacky groupoids.

\begin{lemma} \label{lemma:equi-2gpd-slie} If $f: Y\to X$ is an equivalence, then
\begin{enumerate}
\item \label{itm:i2s} the groupoid $H_1 \rra H_0$ constructed from $Y$ satisfies
$H_1 \cong G_1 \times_{G_0\times_M G_0} H_0\times H_0$ with the
pull-back groupoid structure (therefore $\cG\cong \cH$);
\item \label{itm:algd-morp} the above map $\cG \cong \cH$ induces an SLie groupoid
isomorphism, that is
\begin{enumerate}
\item \label{itm:ii2s}there are a 2-morphism $a: m_{\cG} \to m_{\cH}$:
$\cG\times_M\cG(\cong\cH \times_M\cH)  \to \cG(\cong \cH)$ and a
2-morphism $b:\be_{\cG}\to \be_{\cH}: M \to \cG \cong \cH$;
\item \label{itm:iii2s} between maps $\cG\times_M\cG\times_M\cG \to \cG$, there is a commutative diagram of 2-morphisms,
\[
\xymatrix{
 m_\cG\circ (m_\cG\times id) \ar[r]^{a_G} \ar[d]^{a} & m_\cG\circ
(id\times m_\cG) \ar[d]^{a} \\  m_\cH\circ ( m_\cH \times id)
\ar[r]^{a_H}  & m_\cH\circ (id\times m_\cH) }
\]
where by abuse of notation $a$ denotes the 2-morphisms generated
by $a$ such as $a\odot (a\times id)$;
\item \label{itm:iv2s} between maps $M\times_{M} \cG \to \cG$ and maps
$\cG\times_M M \to \cG$ there are commutative diagrams of
2-morphisms,
\[
\xymatrix{ m_\cG\circ(\be_\cG \times id) \ar[r]^-{b_l^G} \ar[d]^{a\odot b} & pr_2 \\
 m_\cH \circ(\be_cH \times id) \ar[ur]^{b_l^H} & {} } \quad \xymatrix{ m_\cG\circ (id\times \be_\cG
 ) \ar[r]^-{b_r^G} \ar[d]^{a\odot b} & pr_1 \\
 m_\cH \circ(id\times \be_\cH) \ar[ur]^{b^H_r} & {} }
 \]
\end{enumerate}
\end{enumerate}
\end{lemma}
\begin{proof}
Since $Y_2\cong \hom(\partial \Delta^2, Y)\times_{\hom(\partial
\Delta^2, X)} X_2$, we have
\[
\begin{split}
H_1=d_2^{-1}(s_0(Y_0))&=d_2^{-1}(s_0(X_0)) \times_{d_1\times d_0,
X_1\times_M X_1 } Y_1\times_M Y_1 \\
&=G_1 \times_{\bt_G\times \bs_G, G_0\times_M G_0} H_0\times_M H_0 \\
&=G_1\times_{\bt_G\times \bs_G, G_0\times G_0} H_0\times H_0,
\end{split}
\]
where $\bt_G, \bs_G$ denote the target and source map of $G_1\rra
G_0$, and the last step follows from the fact that $(\bt_G \times
\bs_G)(G_1) \subset G_0\times_M G_0$ and $f$ preserves simplicial
structures. The multiplication on $H_1$ (respectively $G_1$) is
given by 3-multiplications on $Y_2$ (respectively $X_2$).
Therefore $H$ has the pull-back groupoid structure since $Y$ is
the pull-back of $X$. So item \eqref{itm:i2s} is proven.

To prove \ref{itm:ii2s}, we translate it to the following groupoid
diagram:
\[
\xymatrix{ & H_1\times_M H_1 \ar[d] \ar@<-1ex>[d] \ar[ld] & & E_{m_\cH}=Y_2 \ar[dll]\ar[drr] & & H_1 \ar[d] \ar@<-1ex>[d] \ar[ld]& \\
G_1\times_M G_1 \ar[d] \ar@<-1ex>[d] & H_0\times_M H_0  \ar[ld] & E_{m_\cG}= X_2 \ar[dll]\ar[drr] & & G_1 \ar[d] \ar@<-1ex>[d] & H_0  \ar[ld] \\
G_0\times_M G_0  & & & & G_0 & }
\]
We need to show that the following compositions of  bibundles are
isomorphic
\begin{equation}\label{eq:2-morp-Em}
\begin{split} &  ( Y_2\times_{G_0}G_1 ) /H_1 (= \big( Y_2 \times_{H_0} H_0 \times_{G_0} G_1 \big) /H_1)
 \\
\overset{a}{\cong}  &H_0\times_M H_0\times_{G_0\times_M G_0} X_2
(=\big( H_0 \times_M H_0 \times_{G_0\times_M G_0}  G_1\times_M G_1
\times_{G_0\times_M G_0} X_2\big) / G_1\times_M G_1) .
\end{split}
\end{equation}
By item \eqref{itm:i2s}, any element in $(Y_2\times_{G_0}G_1)/H_1$
can be written as $[(\eta, 1)]$ with $\eta\in Y_2$, and we
construct $a$ by $[(\eta, 1)] \mapsto (d_2(\eta), d_0(\eta),
f_2(\eta))$. First of all, we need to show that $a$ is
well-defined. For this we only have to notice the identification
$\eta =(\bareta, h_0, h_1, h_2)$ with $\bareta=f_2(\eta)\in X_2$
and $h_i=d_i(\eta)$ since $f_2$ preserves degeneracy maps. Also
$\gamma \in H_1$ can be written as $\gamma=( \bgamma,h_1, h'_1)$
with $\bgamma=f_2(\gamma) \in G_1$, then the $H_1$ action on $Y_2$
is induced in the following way,
\[ (\bareta, h_0, h_1, h_2) \cdot(
\bgamma,h_1, h'_1) = (\bareta \cdot \bgamma, h_0, h'_1, h_2) \]
where $h_i=d_i(\eta)$. Hence if $(\eta',1)=(\eta,1)\cdot(
\bgamma,h_1, h'_1)$, then $\bgamma=1$ and
$a([(\eta',1)])=(h_2,h_0,\bareta)=a([(\eta,1)]$. Given $(h_2, h_0,
\bareta) \in H_0\times_M H_0\times_{G_0\times_M G_0} X_2$, take
any $h_1$ such that $f_1(h_1)=d_1(\bareta)$, then $(h_0, h_1, h_2)
\in \hom(\partial\Delta^2, Y)$. Thus we construct $a^{-1}$ by
$(h_2, h_0, \bareta)\mapsto [((\bareta, h_0, h_1, h_2), 1)]$. By
the action of $H_1$ it is easy to see that $a^{-1}$ is also
well-defined. For the 2-morphism $b$ the proof is much easier
since in this case all the H.S. morphisms are strict morphisms of
groupoids. So we only have to notice the commutative diagram
\[
\xymatrix{M\ar[r] \ar[d] & G_0 \\
H_0 \ar[ru]}
\]

Recall that the 3-multiplications on $Y_2$ are induced from those
of $X_2$ in the following form:
\[ m_0((\bareta_1, h_2, h_5, h_3), (\bareta_2, h_4, h_5, h_0), (\bareta_3, h_1, h_3, h_0))= (m_0(\bareta_1, \bareta_2, \bareta_3), h_2, h_4, h_1), \]
and similarly for other $m$'s.
\[ \xymatrix{ &  & 0  & & \\
1 \ar[urr]^{h_0} & & & & 3 \ar[llll]^{h_4} \ar[llu]^{h_5} \ar[llld]^{h_2} \\
&  2 \ar[ruu]^{h_3} \ar[lu]^{h_1} & & & } \] We also translate
\ref{itm:iii2s} to diagrams of groupoids, where we denote
$n$-product $\square_i\times_M \dots \times_M \square_i$ by
$\square_i^n$,
\[ \small{\xymatrix{ & H_1^3 \ar[d] \ar@<-1ex>[d] \ar[ld] & & Y_2\mathop\times\limits_M H_1 \ar[dll]\ar[drr] & & H_1^2 \ar[d] \ar@<-1ex>[d] \ar[ld]& & Y_2 \ar[dll]\ar[drr]  & & H_1\ar[d] \ar@<-1ex>[d] \\
G_1^3 \ar[d] \ar@<-1ex>[d] & H_0^3  \ar[ld] &  X_2\mathop\times\limits_M G_1 \ar[dll]\ar[drr] & & G_1^2 \ar[d] \ar@<-1ex>[d] & H_0^2  \ar[ld] & X_2 \ar[dll]\ar[drr] & & G_1\ar[d] \ar@<-1ex>[d] & H_0 \\
G_0^3  & & &{H_1\mathop\times\limits_M Y_2}\ar@{.>}[llu]\ar@{.>}[rru]  & G_0^2 & & & & G_0 \\
& & {G_1\mathop\times\limits_M X_2} \ar@{.>}[llu]\ar@{.>}[rru] & &
& & & & & &
 }}\]
Then we have a diagram of 2-morphisms between composed bibundles
\[
\xymatrix{
{\big( (Y_2\times_M H_1 \times_{H_0^2} Y_2 / H_1^2)\times_{G_0}G_1 \big) / H_1} \ar[d]^{a} \ar[r]^{a_H} & {\big( (H_1\times_M Y_2 \times_{G_0^2} Y_2  \big) / H_1^2)\times_{G_0} G_1/ H_1} \ar[d]^{a} \\ \big( H_0^3\times_{G_0^3}(X_2 \times_{M}G_1 \times_{G_0^2} X_2)  \big) /G_1^2
\ar[r]^{a_G} &
\big( (H_0^3\times_{G_0^3}( G_1\times_M X_2 \times_{G_0^2} X_2 ) \big) /G_1^2
}
\]which is
\[ \xymatrix{[(\eta_3, 1, \eta_2, 1)] \ar@{|->}[r]^{a_H} \ar@{|->}[d]^{a} & [(1, \eta_0, \eta_2, 1)] \ar@{|->}[d]^{a} \\
[(h_0, h_1, h_2, \bareta_3, 1, \bareta_1)]  \ar@{|.>}[r]^{a_G ?} &
[(h_0, h_1, h_2, 1, \bareta_0, \bareta_2)]}\] where
$\bareta_i=f_2(\eta_i)$. Since $a_G([\bareta_3, 1,
\bareta_1])=[(1, m_0(\bareta_1, \bareta_2, \bareta_3),
\bareta_2)]$, $f_2$ preserves the 3-multiplications if and only if
$m_0(\bareta_1, \bareta_2, \bareta_3)=\bareta_0$, i.e. if and only
if the above diagram commutes. So \ref{itm:iii2s} is also proven.

Translating the diagrams in item \eqref{itm:iv2s} in groupoid
language, we have
\begin{equation}\label{diag:a-bl-br}
\xymatrix@C=2cm{ \big( J_l^{-1}(M\times_M H_0) \times_{G_0} G_1
\big) / H_1 \ar[d]^{b^H_r}_{b^H_l} \ar[r]^{\text{restriction of
$a$}} &
M\times_M H_0 \times_{M\times_M G_0} J_l^{-1}(M\times_M G_0) \ar[d]^{b^G_r}_{b^G_l} \\
(H_1\times_{G_0}G_1)/H_1 \ar[r]^{f_2} & M\times_M H_0
\times_{M\times_M G_0} G_1}
\end{equation}
where $J_l$ denotes the left moment map of $X_2$ or $Y_2$ to
$G_0\times_M G_0$ or $H_0\times_M H_0$. Then the left diagram of
\eqref{itm:iv2s} is trivial since $b^{H,G}_l=id$ by the
construction in Section \ref{sec:2-slie}. For the right diagram,
we need to show that these two maps are the same: $[(\eta, 1)]
\overset{a}{\mapsto} (h_2, s_0(x), \bareta)
\overset{b^G_r}{\mapsto}(h_2, s_0(x), b_r^G(\bareta))$ and
$[(\eta, 1)] \overset{b_r^H}{\mapsto} [(b_r^H(\eta), 1)]
\overset{f_2}{\mapsto} (h_2, s_0(x), f_2(b_r^H(\eta)))$, where
$x=d_1(h_2)$. That is, we need to show $b_r^G
(f_2(\eta))=f_2(b_r^H(\eta))$. Since $b_r=\varphi^{-1}$ is
constructed by $m$'s as in Section \ref{sec:2-slie}, $f_2$
commutes with $b_r$'s. So we proved item \eqref{itm:iv2s}.
\end{proof}

To establish the inverse argument, we fix again the notation:
$\cG\rra M$ is an SLie groupoid; $G_1\rra G_0$ and $H_1 \rra H_0$
are two Lie groupoids presenting $\cG$; $X$ and $Y$ are the
2-groupoids corresponding to $G$ and $H$ as constructed in Section
\ref{sec:slie-2}.

\begin{lemma}If there is a Lie groupoid equivalence $H \to G$,
then there is a Lie 2-groupoid equivalence $Y\to X$.
\end{lemma}
\begin{proof}
Since both $H$ and $G$ present $\cG$ which is an SLie groupoid
over $M$, we are in the situation described in item 2 in Lemma
\ref{defi:1t-m-equi-2gpd}, that is we have a 2-morphism $a$
satisfying various commutative diagram as in item
\eqref{itm:ii2s}, \eqref{itm:iii2s}, \eqref{itm:iv2s}. We take
$f_0$ to be the isomorphism $M\cong M$, $f_1$ the map $H_0\thra
G_0$, $f_2: Y_2\to X_2$ the map $\eta \mapsto [(\eta, 1)]
\overset{a}{\mapsto} (h_2, h_0, \bareta) \mapsto \bareta$. Since
$f_2$ is made up by composition of smooth maps, it is a smooth
map. Since $d_2\times d_0$ is the moment map and $a$ preserves the
moment map, we have $h_i=d_i(\eta)$ for $i=0,2$. It is clear that
$f_0$ and $f_1$ preserve the structure maps since they preserve
$\be$, $\bbs$, $\bbt$ of $\cG\rra M$. It is also clear that $d_i
f_2(\eta)=f_1(h_i)$ for $i=2,0$ since $(h_2, h_0,
\bareta=f_2(\eta)) \in H_0^2\times_{G_0^2} X_2$. Since $a$
preserves moment maps, $f_1(d_1(\eta))=J_r([(\eta,
1)])=J_r(h_2,h_0,\bareta)=d_1 (\bareta)$, where $J_r$ is the
moment map to $G_0$ of the corresponding bibundles. Hence $f_2$
preserves the degeneracy maps.

For the face maps $s_0, s_1: \square_1\to \square_2$, we again
look at the commutative diagram \eqref{diag:a-bl-br}. Recalling
$s_0(h)=b_l^{-1} \be_H (h)$, we have
\[ f_2(s_0(h))= pr_X a([(b_l^{-1} \be_H(h), 1)]) = b_l^{-1} \be_\cG
f_1(h) = s_0 f_1(h), \]where $pr_X$ is the projection
$H_0^2\times_{G_0^2} X_2 \to X_2$ and similarly for $s_1$ using
the diagram for $b_r$. Hence $f_2$ preserves the face maps.

The fact that $f_2$ preserves the 3-multiplications follows from
the proof of item \eqref{itm:iii2s} of Lemma
\ref{lemma:equi-2gpd-slie}.

Then the induced map $\phi: Y_2\to \hom(\partial \Delta^2, Y)
\times_{\hom(\partial \Delta^2, X)} X_2$ is $\eta \mapsto
[(\eta,1)]\overset{a}{\mapsto} (h_2, h_0, \bareta) \mapsto (h_0,
h_1, h_2, \bareta)$, where $\bareta=f_2(\eta)$ and
$h_i=d_i(\eta)$. As a composition of smooth maps $\phi$ is smooth.
Moreover $\phi$ is injective since $a$ is injective. For any
$(h_0, h_1, h_2,\bareta)\in \hom(\partial\Delta^2, Y)
\times_{\hom(\partial\Delta^2, X)} Y_2$, we have $(h_0, h_2,
\bareta)\in H_0\times_M H_0\times_{G_0\times_M G_0} X_2$. Then
there is an $\eta$ such that $[(\eta, 1)]=a^{-1}(h_0, h_2,
\bareta)$. Thus $\phi (\eta)=(h_0, d_1(\eta), h_2, \bareta) $,
which implies that $f_1 (d_1(\eta))=d_1(\bareta)=f_1(h_1)$.
Therefore $(1, d_1(\eta), h_1)\in H_1$ and $d_i(\eta \cdot (1,
d_1(\eta),h_1))=h_i$, $i=0,1,2$ since $d_1$ is the moment map to
$H_0$ of the bibundle $Y_2$. So $\phi(\eta \cdot (1, d_1(\eta),
h_1))= (\bareta, h_0, h_1, h_2)$ which shows the surjectivity.
Therefore $\phi$ is an isomorphism.
\end{proof}

The theorem is now proven, since we only have to consider the case
that when (1-) Morita equivalence is given by Lie (2-) groupoid
morphisms.

\section{Comments on the inverse map}\label{sec:inverse}
In this section, we prove that the axioms involving the inverse
map in the definition of SLie groupoid are guaranteed by the
axioms of multiplication and the identity. There is similar
treatment of the antipode in hopfish algebras \cite{twz}. In fact
SLie groups are a geometric version of hopfish algebras.

Let $\cG\rra M$ be an SLie groupoid, and $G:=G_1
\underset{\bt_G}{\overset{\bs_G}{ \rra }}G_0$ a good groupoid
presentation of $\cG$ as we described in Section
\ref{sec:embedding}. 
So there is a map $e: M\to G_0$ presenting $\bar{e}$.  We look at
the following diagram,
\[ \cG \times_{M} \cG \overset{m}{\lra}\cG \overset{\be}{\lla}M, \]
and its corresponding groupoid picture,
\begin{equation}\label{eq:inv-construct}
\xymatrix{ G_1\times_M G_1 \ar[dd]\ar@<-1ex>[dd]&  &G_1\ar[dd]\ar@<-1ex>[dd]& &M\ar[dd]\ar@<-1ex>[dd]\\
& E_m\ar[dr]^{J_r}\ar[dl]_{J_l}& & {E_{\be}} \ar[dr]^{J_r}\ar[dl]_{J_l}& &\\
G_0\times_M G_0&  &G_0& &M\\
}
\end{equation}
where $E_m$ and $E_{\be} =G_1\times_{\bt_G, G_0, e}M$ 
are bibundles presenting the multiplication $m$ and identity $\be$
of $\cG$ respectively. We can form a $G\times_M G$ left module
$E_m \times_{G_0} E_{\be} /G$.  Examining the $G$ action on
$E_{\be}$, we see that the geometric quotient (corresponding to
$\hom_\cA (\bepsilon, \bDelta)$ in the case of hopfish algebra),
\[ (E_m \times_{G_0} E_{\be})/G = E_m\times_{J_r, G_0, e} M, \]
is still a manifold with a left $G_1\times_M G_1$ action (which
might not be free and proper). Therefore, we can view it as a left
$G$ module with the left action of the first copy of $G$ and a
right $G^{op}$ module with the left action of the second copy of
$G$. Here $G^{op}$ is $G$ with the opposite groupoid structure.

\begin{lemma} \label{proof-morita} The bibundle above is a Morita
 equivalence
 from $G$ to $G^{op}$ with moment maps $pr_1\circ J_l$ and $pr_2\circ J_r$.
\end{lemma}
\begin{proof}
 Let $f_1: \cG \times_{M} \cG \to \cG \times_{M} \cG$ be given by $(g_1, g_2)\mapsto (g_1 \cdot g_2, g_2)$, i.e. $f_1 = m\times pr_2$; let $f_2: \cG \times_M \cG \to \cG \times_M \cG$ be given by $(g_1, g_2)\mapsto (g_1 \cdot g_2^{-1}, g_2)$. Since we have
\[ (g_1, g_2) \overset{f_1}{\mapsto} (g_1 g_2, g_2)\overset{f_2}\mapsto ((g_1 g_2)g_2^{-1}, g_2))\sim (g_1, g_2), \]
and
\[ (g_1, g_2) \overset{f_2}{\mapsto}(g_1 g_2^{-1}, g_2)\overset{f_1}{\mapsto}((g_1 g_2^{-1})g_2 , g_2))\sim (g_1, g_2), \]
we have that $f_1 \circ f_2$ and $f_2\circ f_1$ are isomorphic to
$id$ via 2-morphisms. Therefore $f_1$ is an isomorphism of stacks.
Therefore $E_m \times_{pr_2 \circ J_l, G_0, \bt_G} G_1$ presenting
$f_1$ is a Morita bibundle from the Lie groupoid $G_1 \times_M
G_1\rightrightarrows G_0 \times_M G_0$ to $G_1\times_M G_1
\rightrightarrows G_0\times_M G_0$. Here the two moment maps are
$J_l$ (of $E_m$) and $J_r \times \bs_G$.  Therefore the left
groupoid action of  $G_1 \times_M G_1$ is principal on the
bibundle $E_m \times_{pr_2 \circ J_l, G_0, \bt_G} G_1$. The left
action of the first copy of $G_1$ on $E_m$ is the left action we
need. Therefore, it is naturally free. To see that the desired
action is also transitive, take $(\eta_1, x)$ and $(\eta_2, x)$ in
$E_m\times_{J_r, G_0, e} M$ which project to the same point $g_0$
in $G_0$ by  $pr_2\circ J_l$. Then $(\eta_1, 1_{g_0})$ and
$(\eta_2, 1_{g_0})$ are both in  $E_m \times_{pr_2 \circ J_l, G_0,
\bt_G} G_1$ and both project to $(e(x), g_0)\in
 G_0 \times_M G_0$ via $J_r \times \bs_G$. Therefore, there is an element
$(\gamma_1, \gamma_2)\in G_1\times_M G_1$ such that $(\gamma_1,
\gamma_2)\cdot (\eta_1, 1_{g_0})= (\eta_2, 1_{g_0})$. Therefore
$\gamma_1 \cdot (\eta_1, x)=(\eta_2, x)$, i.e. the action is
transitive fibre-wise.

The proof of the principality of the  $G^{op}$ action is similar
(one considers $G_1\times_{G_0}E_m$).

\end{proof}
\begin{remark}
Another fibre product $E_m \times_{pr_2\circ J_l, G_0, e}M $ is
isomorphic to $ G_1$ trivially via $b_r$. But the morphisms we use
to construct the fibre product are different.
\end{remark}

Notice that using the inverse operation, a $G^{op}$ module is also
a $G$ module. In other words, the above lemma says that
$E_m\times_{J_r, G_0, e} M$ is a Morita bibundle between $G$ and
$G$ where the right $G$ action is via the left action of the
second of copy of $G\times_M G$ composed with the inverse. With
this viewpoint, we have a stronger statement:
\begin{lemma} \label{isom}
As Morita bibundles from $G$ to $G$, $E_m\times_{J_r, G_0, e} M$
and $E_i$ are isomorphic.
\end{lemma}
This lemma is not at all surprising, if we use the view point of
Lie 2-groupoids: $E_m\times_{J_r, G_0, e} M$ should be imagined as
degenerate triangles (on the left)
\begin{equation} \label{diag:inverse}
\begin{xy}
*\xybox{(0,0);<3mm,0mm>:<0mm,3mm>::
  ,0
  ,{\xylattice{-5}{0}{-4}{0}}}="S",
  {(-10,-10)*{\bullet}}, {(-10, -12)*{_1}},
     {(12,-7)*{\bullet}}, {(14, -7)*{^{0}}},
     {(10, -10)*{\bullet}}, {(10, -12)*{^{2}}},
     {(-10, -10) \ar@{->} (12,-7)}, 
     { (10, -10) \ar@{->} (-10,-10)}, { (10, -10) \ar@{->} (12,-7)} , 
\end{xy} \quad
\begin{xy}
*\xybox{(0,0);<3mm,0mm>:<0mm,3mm>::
  ,0
  ,{\xylattice{-5}{0}{-4}{0}}}="S",
  {(-10,-10)*{\bullet}}, {(-10, -12)*{_1}},
     {(12,-7)*{\bullet}}, {(14, -7)*{^{0}}},
     {(10, -10)*{\bullet}}, {(10, -12)*{^{2}}},
     {(10, -4)*{\bullet}}, {(14,-4)*{^{0'<0}}},
     {(-10, -10) \ar@{->} (12,-7)}, 
     { (10, -10) \ar@{->} (-10,-10)}, { (10, -10) \ar@{->} (12,-7)} , 
      {(10, -4)\ar@{<-}(-10, -10)}, 
      {(10, -4)\ar@{<-}(12, -7)}, 
      {(6, -7)*{^{\gamma_0}}}, {(7, -9)*{^{\eta_1}}},{(0, -10)*{^{\eta_2}}}
\end{xy}
\end{equation}

\begin{proof}
We know from the property of $E_i$ that $g\cdot g^{-1} \sim 1$,
i.e. $((G_1\times_{\bt_G,G_0,J_l}E_i)\times_{\bs_G \times J_r,
G_0\times_M G_0}E_m) /G_1\times_M G_1$ is 2-isomorphic to the
module $G_0\times_{e\circ \bt, G_0, \bt_G}G_1$ presenting the map
$e\circ \bbt: \cG\to M \to \cG$.  We will first show that $E_m
\times_{G_0}M$ also has this property.

Let $[(1, \eta_1, \eta_2)]\in (((G_1\times_{G_0}(E_m\times_{J_r,
G_0, e} M))\times_{G_0\times_M G_0}E_m)/G_1\times_M G_1$ (See
\eqref{diag:inverse}). We can suppose the first element to be 1
because the $G_1 \times_M G_1$ action on $G_1$ is simply right
multiplication by the first factor. Moreover, examining the
morphisms in the fibre product, we have $J_l(\eta_1)=J_l(\eta_2)$,
and the right $G_1\times_M G_1$ action is
\begin{equation}\label{eq:r-gg-act} (\gamma_0, \eta_1, \eta_2)\cdot
(\gamma_1, \gamma_2)= (\gamma_0 \cdot \gamma_1, (1,
\gamma_2^{-1})\cdot \eta_1, (\gamma_1, \gamma_2)^{-1} \cdot
\eta_2).\end{equation} Since the right action of $G_1$ on $E_m$ is
principal (now viewing $E_m$ as a bibundle from $G\times_M G$ to
$G$), and $J_l(\eta_1)=J_l(\eta_2)$, there exists a unique element
$\gamma\in G_1$ such that $\eta_1\cdot \gamma=\eta_2$. Suppose
$[(1, \eta_1, \eta_2)]=[(1,\teta_1, \teta_2)]$, then by
\eqref{eq:r-gg-act}, there exists $\gamma_2 \in G_1$ such that
$\eta_1= (1, \gamma_2^{-1})\cdot \teta_1 $ and $\eta_2=(1,
\gamma_2^{-1})\cdot \teta_2$. Since the right action and left
action on a bibundle commute, $\gamma$ is independent of the
choice of representative of $[(1, \eta_1, \eta_2)]$. Therefore, we
have a well-defined map
\[ \phi:((G_1\times_{G_0}(E_m\times_{J_r, G_0, e} M))\times_{G_0\times_M G_0}E_m)/G_1\times_M G_1 \to   G_0\times_{e\circ \bt, G_0, \bt_G}G_1,\]
by
\[ [(1_g, \eta_1, \eta_2)]\mapsto (g, \gamma). \]
Here $(g, \gamma)$ is an element of $ G_0\times_{e\circ \bt, G_0,
\bt_G}G_1$ since $\bt_G(\gamma)= J_r (\eta_1)$. Moreover, since
smoothness is a local property and  all the moment maps are smooth
and the groupoid action is principal, $\phi$ is naturally smooth.
$\phi$ is an isomorphism following from  principality of the right
$G_1$ action on $E_m$. Moreover, it is not hard to check that
$\phi$ is equivariant and  commutes with the moment maps of the
bibundles. Therefore,
\[ ((G_1\times_{G_0}(E_m\times_{J_r, G_0, e} M))\times_{G_0\times_M G_0}E_m)/G_1\times_M G_1 \cong  G_0\times_{e\circ \bt, G_0, \bt_G}G_1 \]as Morita bibundles.
One proceeds similarly to prove the other symmetric isomorphism
corresponding to $g^{-1}\cdot g \sim 1$.

Let $\varphi$ be the composed isomorphism
\begin{equation}\label{compose}
((G_1\times_{G_0}(E_m\times_{J_r, G_0, e} M))\times_{G_0\times_M
G_0}E_m)/G_1\times_M G_1 \to
((G_1\times_{G_0}E_i)\times_{G_0\times_M G_0}E_m) /G_1\times_M
G_1.
\end{equation}
Suppose $\varphi([(1_g, \eta_1, \eta_2)])=([(1_g, \teta_1,
\teta_2)])$ (we can still assume that the first component is 1
because the $G_1\times_M G_1$ action on both side is the right
multiplication by the first copy; we can assume that they are 1 at
the same point because $\varphi$ commutes with the moment maps on
the left leg). Examining the morphisms inside the fibre products,
we have
\[ pr_1 \circ J_l (\eta_2)=\bt_G(1_g)=pr_1\circ J_l(\teta_2)=g. \]
Since $\varphi$ commutes with the moment maps on the right leg, we
have
\[ J_r(\eta_2)=J_r(\teta_2).\]
Similarly to the proof of Lemma \ref{proof-morita}, we can show
that $G_1\times_{\bs_G, G_0, pr_1\circ J_l} E_m$ is a Morita
bibundle from $G \times_M G$ to $G \times_M G$. Then $(1_g,
\eta_2)$ and $(1_g, \teta_2)$ are both in $G_1\times_{\bs_G, G_0,
pr_1\circ J_l} E_m$ and their images under the right moment map
$\bs_G\times J_r$ are both $(g, J_r(\eta_2))$. By principality of
this left $G_1\times_M G_1$ action, there is a unique $(\gamma_1,
\gamma_2)\in G_1\times_M G_1$ such that
\[ (\gamma_1, \gamma_2)\cdot (1, \eta_2)=(1, \teta_2). \]
Therefore $\gamma_1=1$ and $(1, \gamma_2) \cdot \eta_2 =\teta_2$.
This left $G_1 \times_M G_1$ action on $E_m$ is exactly the left
$G_1\times_M G_1$ action on the second copy of $E_m$ in
\eqref{compose}. Using this $\gamma_2$, we have
\[ (1, \teta_1, \teta_2) \cdot (1, \gamma_2) =(1, \gamma_2^{-1})\cdot (1,\teta_1,\teta_2)= (1, \eta_1', \eta_2).\]
Therefore the isomorphism
\[\varphi: \; [(1_g, \eta_1, \eta_2)]\mapsto [(1_g, \eta'_1, \eta_2)]\]
induces a map $\psi: E_m\times_{G_0} M \to E_i$ by $ \eta_1
\mapsto \eta_1'$. It's routine to check $\psi$ is an isomorphism
of Morita bibundles.
\end{proof}

As in this Lemma, we have seen that the 2-identities satisfied by
$E_i$ are actually naturally satisfied by $E_m \times_{J_r, G_0,
e}M$. Notice that for the first part of the proof, we didn't use
any information involving the inverse map.  Our conclusion is that
the inverse map represented by $E_i$ can be replaced by $E_m
\times_{J_r, G_0, e}M$ without any further conditions (not even on
the 2-morphisms). In fact, the natural 2-morphisms coming along
with bibundle $E_m \times_{J_r, G_0, e}M$ naturally go well with
other 2-morphisms $a$'s and $b$'s.

\begin{prop} \label{prop:inverse}
An SLie groupoid $\cG$ can be equally defined by all the axioms
except the axioms involving inverses and the additional condition
that $E_m \times_{J_r, G_0, e}M$ is a Morita bibundle from $G$ to
$G^{op}$ for some good presentation $G$ of $\cG$.
\end{prop}
\begin{proof} It is clear from the Lemma \ref{isom} that the existence of the
inverse map guarantees that the bibundle  $E_m \times_{J_r, G_0,
e}M$ is a Morita bibundle from $G$ to $G^{op}$ for a good
presentation $G$ of $\cG$.

On the other hand, if  $E_m \times_{J_r, G_0, e}M$ is a Morita
bibundle from $G$ to $G^{op}$ for some presentation $G$ of $\cG$,
then we construct the inverse map $i: \cG \to \cG$ by this
bibundle. Because of the nice properties of  $E_m \times_{J_r,
G_0, e}M$ that we have proven in the first half of Lemma
\ref{isom}, this newly defined inverse map satisfies all the
axioms that the inverse map satisfies. \end{proof}

\begin{remark}
This theorem holds also for W-groupoids and the proof is similar.
Then one can see that the new definition of SLie group modulo
2-morphisms is analogous to the definition of hopfish algebra.
\end{remark}

Sometimes the inverse of an SLie groupoid is given by a specific
groupoid isomorphism $i: G \to G$ on some presentation (for
example $\cG(A)$ and $\cH(A)$ in \cite{tz} and (quasi-)Hopf
algebras as the algebra counter-part).

\begin{lemma}
The inverse map of an SLie groupoid $\cG$ is given by a groupoid
isomorphism $i: G\to G$ for some presentation $G$ if and only if
on this presentation $E_m \times_{J_r, G_0, e}M$ is a trivial
right $G$ principal bundle over $G_0$.
\end{lemma}
\begin{proof} It follows from Lemma \ref{isom} and the fact that the inverse is given by a morphism $i: G\to G$ if only if
the bibundle $E_i$ is trivial.
\end{proof}

\bibliographystyle{alpha}
\bibliography{bibz}

\end{document}